\documentclass{article}
\usepackage[utf8]{inputenc}
\usepackage{amsmath}
\usepackage{amsfonts}
\usepackage{amssymb}
\usepackage{latexsym}
\usepackage[titletoc]{appendix}
\newtheorem{te}{Theorem}[section]
\newtheorem{lemma}{Lemma}[section]
\newtheorem{coro}{Corollary}[section]
\newtheorem{rem}{Remark}[section]
\newcommand{\cc}{\mathbb{C}}
\newcommand{\rr}{\mathbb{R}}
\newcommand{\va}{\varphi}
\newcommand{\al}{\alpha}
\newcommand{\pp}{\partial}
\newcommand{\dis}{\displaystyle}
\newcommand{\la}{\lambda}
\usepackage{hyperref}
\hypersetup{
  hidelinks,
  urlcolor=black,
}
\begin{document}
\title{\bf  Non-self adjoint Sturm-Liouville problem with spectral  and physical parameters in boundary conditions}
\author{Rodrigo Meneses Pacheco\\
\\
Escuela de Ingenier\'ia Civil, \\
Facultad de Ingenier\'ia, Universidad de Valpara\'iso\\
 Avda.
Errazuriz 1834, Valpara\'iso, Chile\\
\href{mailto:rodrigo.meneses@uv.cl}{rodrigo.meneses@uv.cl}\\
\\
Oscar Orellana\\
\\
Departamento de Matem\'aticas\\
Universidad T\'ecnica Federico santa Mar\'ia\\
Avenida Espa\~na 1680, Valpara\'iso, Chile\\
\href{mailto:oscar.orellana@usm.cl}{oscar.orellana@usm.cl}
}
\maketitle
\begin{abstract}
We present a complete description on the spectrum and eigenfunctions of the following two point boundary value problem
\begin{equation}
\label{pabstract}
\left\lbrace\begin{array}{rcl}
 \dis{(p(x)f')'-(q(x)- \lambda r(x))f}&=&0\qquad 0<x<L \\
                                          f'(0)&=&(\alpha_{1} \lambda +\alpha_{2})f(0)\\
                                          f'(L)&=&(\beta_{1}\lambda -\beta_{2})f(L)
            \end{array}\right.
\end{equation}
where $\la$ and $\al_{i},\ \beta_{i}$ are spectral and physical parameters.
Our survey is focused mainly in the case $\alpha_{1}>0$ and $\beta_{1}<0$, where neither self adjoint operator theorems on Hilbert spaces nor Sturm's comparison results can be used directly.
We describe the spectrum and the oscillatory results of the eigenfunctions from a geometrical approach, using a function related to the Pr\"ufer angle.
The proofs of the asymptotic results of the eigenvalues and separation theorem of the eigenfunctions are developed through classical second order differential equation tools.
Finally, the results on the spectrum of \eqref{pabstract} are used for the study of the linear instability of a simple model for the fingering phenomenon on the flooding oil recovery process.
\end{abstract}

\section{Introduction}
It is well known that the results on boundary valued problems have a direct and relevant implication on many models in applied mathematics, as can be seen for the vast literature on the subject \cite{chandrasekhar1961hydrodynamic,amrein2005sturm,drazinhydrodynamic}.
On many occasions, the boundary conditions are described as functions of the spectral parameter and other physical parameters.
Problems of this type are found in mechanic models \cite{chandrasekhar1961hydrodynamic,rayleigh1896theory,timoshenko1974vibration,budak1988collection}.
On elastic models of rod and string they can usually be described through self adjoint Sturm-Liouville problem (SLP) \cite{langer1947fourier} while that the spectral problems arising in hydrodynamics are usually  non-self adjoint \cite{chandrasekhar1954characteristic,mennicken2003non}.
An extensive classical bibliography for various physical applications can be found in \cite{amara1999sturm,fulton1977two,walter1973regular}.

In this article we study the following SLP with spectral parameter in the boundary conditions
\begin{equation}
\label{problema}
\left\lbrace\begin{array}{rcl}
 \dis{(p(x)f')'-(q(x)- \lambda r(x))f}&=&0\qquad 0<x<L \\
                                          f'(0)&=&(\alpha_{1} \lambda +\alpha_{2})f(0)\\
                                          f'(L)&=&(\beta_{1}\lambda -\beta_{2})f(L)
            \end{array}\right.
\end{equation}
where $\la$ is a spectral parameter, and  $\al_{i}$ and $\beta_{i}$ can be considered as physical parameters.
Our interest in the study of problem \eqref{problema} arises from the linear stability analysis of a three-layer Hele-Shaw cell model for the  study of a secondary oil recovery process.
This problem is presented in \cite{gorell1983theory} and has been studied in several articles on the subject, see \cite{carasso1998optimal} and references therein.
The model problem is presented in section \ref{applications}, where problem \eqref{problema} corresponds to a synthesis of the physical laws considered.
In the stability model, $\al_{i}$ and $\beta_{k}$ are defined as functions of the wave numbers of the perturbative wave, therefore the signs of those parameters are related to different wave number ranges.

Referring the SLP with spectral parameter on the boundary conditions, the way in which the spectral parameter relates with the physical parameters introduces different types of complications on the tools used for the analysis.

Concerning the tools for the analysis of \eqref{problema}, the sign of parameters $\al_{1}$ and $\beta_{1}$ plays an essential role.
For the case where $\al_{1}<0$ and $\beta_{1}>0$, the description of the spectrum, asymptotic results, oscillatory results on the eigenfunctions and eigenfunction expansion results have been widely developed in several articles (see \cite{churchill1942expansions,cohen1966integral}) and can be learned from classical texts, like Ince \cite{ince1962ordinary} and Reid \cite{reid1971ordinary}, among others.
In this case, other than Sturm's comparison results, it's also possible to use results on selfadjoint operators in Hilbert spaces \cite{fulton1977two}.
The identification and characterization of the spectrum of \eqref{problema} for the case where $\al_{1}\beta_{1}\geq0$ is not as straightforward, since the geometrical approach through the transformation of Pr\"ufer presents new difficulties.
In \cite[Theorem 2.1]{binding2004transformation} the authors describe in detail the behavior of the spectrum of \eqref{problema} with $\al_{1}=0$ and $\beta_{1}<0$ via multiple Crum-Darboux transformations, obtaining an associated 'almost' isospectral regular SLP.
For the case where $\al_{1}=\al_{2}=\beta_{2}=0$, with $\beta_{1}<0$, in the article \cite{amara1999sturm}, the authors address the study of the spectrum (Theorem 2) and and asymptotic results of the eigenvalues and eigenfunctions (Theorem 1) using self adjoint operators in Pontryagin space $\Pi_{1}$.
We note that in the general case, where $\al_{1}\beta_{1}\geq0$, it's possible to apply Sturm's oscillatory results for one of the boundary conditions, which eases the geometric analyses commonly used for the study of the spectrum.
General information on SLP with spectral parameters can be found in \cite{behrndt2006sturm, binding2004transformation, binding2003hierarchy, binding2006sturm, binding1994sturm, churchill1942expansions, ince1962ordinary, reid1971ordinary}.

The main objective of this work is to obtain results for \eqref{problema} when $\al_{1}>0$ and $\beta_{1}<0$, namely on the descriptions and behaviors of the eigenvalues and eigenfunctions, similar to existing results on regular SLP.
In our analysis, we mainly use classical tools on second order ordinary differential equations (ODE), adapting several characteristic results on regular SLP, for example  theorems of Sturm (see Hartman \cite{hartman1964ordinary} 11.3).

For the description of the spectrum and the oscillatory results of the eigenfunctions, we used a geometrical approach through a function related  to the Pr\"ufer angle.
For the results of separation, we make an analysis of the monotone behavior of the zeros of the solution functions, obtaining this result through the implicit function theorem.
Concerning the functions $p(x),\ q(x),\ r(x)$, we assumed them to be positive and sufficiently regular in $[0,L]$, so they allow the use of classical results on regularity of second order ODEs solutions with respect to the model parameters (See Peano's theorem in Hartman \cite{hartman1964ordinary}).
When studying the model in Section \ref{applications}, we also assume that the function $p(x)$ is strictly increasing and that constants $\al_{2}$ and $\beta_{2}$ are positive.

Our main result read as follows:
\begin{te}
\label{principal}
Under the conditions  $\alpha_{1}>0$ and $\beta_{1}<0$, the  spectrum of \eqref{problema} consists of an unbounded sequence of real eigenvalues
\begin{equation}
\label{orden}
\la_{-1}<\la_{-0}<0<\la_{0}<\la_{1}<\la_{2}<\dots\nearrow \infty
\end{equation}
and the corresponding eigenfunction $f(x;\la_{l})$ has exactly $\vert l\vert$ zeros in $]0,L[$.
\end{te}

We prove this using classical geometric tools related with the Pr\"ufer angle.
The technique consists in considering $f(x;\la)$ a solution of the ODE in \eqref{problema}, satisfying the boundary condition at $x=0$. The eigenvalues of \eqref{problema} are studied from the equation

\begin{equation}
\label{laecuacion}
p(L)\frac{f'(L;\la)}{\la f(L;\la)}=\beta_{1}-\frac{\beta_{2}}{\la},
\end{equation}
for those $\la$ where the expression are well defined.
In the course of the proof, the monotony of the function in \eqref{laecuacion} is a fundamental tool.
For this, it's necessary to characterize the branches where the function $h_{1}(\la)=p(L)\frac{f'(L;\la)}{\la f(L;\la)}$ is well defined.
Knowing the behavior of function $h_{1}(\la)$, the solutions of equation \eqref{laecuacion} are studied as the graph intersection of functions $h_{1}(\la)$ and $h_{2}(\la)=\beta_{1}-\frac{\beta_{2}}{\la}$.
When considering the case for $\alpha_{1}>0$ and $\beta_{1}<0$ we can't apply Sturm's first comparison theorem directly, therefore we can't use as corollary Sturm's separation theorem.
The monotony results of $h_{1}(\la)$ on each branch where defined are obtained from the analysis of an ODE, where $h_{1}(\la)$ is a solution.

From the consideration $\al_{1}>0$ and $\beta_{1}<0$, we can't use Sturm's separation theorem.
Next, we present our second main result, which addresses this problem.
\begin{te}[Separation Theorem]
\label{principal2}
For $0<\underline{\la}<\overline{\la}$, consider $\underline{f}(x)=f(x;\underline{\la})$
  and $\overline{f}(x)=f(x;\overline{\la})$.
Let $<x_{1}<x_{2}<L$ two consecutive zeros of $\underline{f}$.
Then $\overline{f}$ has at least one zero on $[x_{1},x_{2}]$.
\end{te}
The proof is made using the behavior of the zeros of function $f(x;\la)$.
For this, we use the implicit function theorem and, therefore, use regularity of coefficients of the ODE in \eqref{problema}, for  the regularity of parameter $\la$.
Using Theorem  \ref{principal2} and the oscillatory results of Sturm, we get the following results on the eigenfunctions directly.
\begin{coro}
\label{separacion}
Let $f_{n}(x)=f(x;\la_{n})$ and $f_{n+1}(x)=f(x;\la_{n+1})$ with $0<\la_{n}<\la_{n+1}$ eigenvalues of \eqref{problema}.
Then between two consecutive zeros of $f_{n+1}$ there is exactly one zero of $f_{n}$.
\end{coro}
We note that for the case where $\la<0$ we have a non-oscillatory ODE, on which other tools, like the maximum principle are used for its analysis.
All the results presented are studied from a geometrical approach, specifically, studying the  behavior of the function $h_{1}(\la)$.

Finally, to understand the qualities of the eigenvalues, we present the following asymptotic result:
\begin{te}
\label{asintotico}
For $\la_{n}$ in \eqref{orden} we have $\sqrt{\la_{n}}=n\pi/L +O(n^{-1})$ as $n\to\infty$
\end{te}
The development of the proof of this theorem is conditioned to the existence of an eigenvalue $\la_{0}$ such that $f(x,\la_{0}) > 0$ in $[0,L]$.
For the proof we use Crum-Darboux and Liouville transformations, obtaining an associated regular SLP where we can relate the spectra of both problems.
The ideas and terminology where taken from article \cite{binding2004transformation} and some of its references.

The present work is ordered as follows: Section \ref{pre1} of preliminaries, presents a list of antecedents used for our main results.
In the subsections we display elemental tools of our analysis.
Such tools correspond (mostly) to adaptations of classical results on regular SLP.
We decided to present a sequence of lemmas for easing the read of the main proofs.
Of the results presented in Section \ref{pre1}, we highlight the behavior of the function $h_{1}(\la)$ defined in \eqref{funciones} and the behavior of the zeros of the functions $f(x;\la)$, solutions for the ODE in \eqref{problema} satisfying the boundary condition for $x=0$.
In Section \ref{demoprincipales1} we develop the proofs of our main results.
In section \ref{demoprincipales2} we use a regular SLP to characterize the asymptotic behavior of the eigenvalues of problem \eqref{problema}.
The results of this work are applied to an hydrodynamics model in Section \ref{applications}.
The problem considered in this section corresponds to a non regular SLP associates to the study of a linear stability problem on a secondary oil recovery process.
The hypotheses on the coefficient functions and parameters in \eqref{problema} correspond to considerations on the studied model.
In this section we develop numeric computations on a particular case of \eqref{problema}.

\section{Preliminaries.}\label{pre1}
In this section we present the main ideas and some results used in the proofs of our theorems.
The aim is to point out some of technical difficulties for the attainment of results for the case where $\al_{1}>0$ and $\beta_{1}<0$.

In order to ease the reading of this section, we indicate some particular cases of our main results separately, specifically, the result of existence of the main eigenvalue $\la_{0}>0$ satisfying $f(x;\la_{0})>0$ in $]0,L[$.
We also note that the existence of such eigenvalue is fundamental for the construction of the proofs for the asymptotic representation of the eigenvalues, presented in Theorem \ref{asintotico}.

We alse present results on the behavior of the zeros of the solutions for the ODE in \eqref{problema} and the first boundary condition (see problem \eqref{main2}).

\subsection{ Geometrical Framework.
Shooting technique  and an auxiliary non-regular SLP        }
We use the approach considered for the description of the spectrum of \eqref{problema} and the oscillatory results of the eigenfunctions.

Under conditions of regularity of the ODE coefficients, there are two functions $f_{1}(x;\lambda)$ and $f_{2}(x;\lambda)$ of class $\mathcal{C}^{2}$, such that they are linearly independent solutions of the ODE in (\ref{pabstract}).
Thus, the general solution of the ODE in (\ref{pabstract}) can be represented as follows
$$ f(x;\lambda)=C_{1}f_{1}(x;\lambda) + C_{2}f_{2}(x;\lambda).$$
Now, applying the first boundary condition of problem (\ref{pabstract}) we get that:
$$ C_{1}f_{1}'(0;\lambda)+C_{2}f_{2}'(0;\lambda)=(\alpha_{1}\lambda + \alpha_{2})(    C_{1}f_{1}(0;\lambda)+C_{2}f_{2}(0;\lambda)         )     $$
and therefore,
$$ C_{2}\left( f_{2}'(0;\lambda)- (\alpha_{1}\lambda + \alpha_{2}) f_{2}(0;\lambda)  \right) = C_{1}(  (\alpha_{1}\lambda + \alpha_{2}) f_{1}(0;\lambda)-f'_{1}(0;\lambda)   ) .
$$
Thus, considering
$$ C_{2}=C_{1}\left(   \frac{(\alpha_{1}\lambda + \alpha_{2}) f_{1}(0;\lambda)-f'_{1}(0;\lambda) }{f_{2}'(0;\lambda)- (\alpha_{1}\lambda + \alpha_{2}) f_{2}(0;\lambda)}          \right)    $$
we obtain the following representation:
$$ f(x;\lambda)= C_{1}\left(   f_{1}(x;\lambda)-\left[  \frac{  f'_{1}(0;\lambda)-(\alpha_{1}\lambda + \alpha_{2}) f_{1}(0;\lambda) }{f_{2}'(0;\lambda)- (\alpha_{1}\lambda + \alpha_{2}) f_{2}(0;\lambda)}     \right]f_{2}(x;\lambda)\right)   $$
where $C_{1}\neq 0$ is an arbitrary constant.

We note that the expression $f'(x;\la)/f(x;\la)$ is completely independent of the value  $C_{1}$.
On the other hand, from the regularity of the coefficient functions $p(x)$, $q(x)$ and $r(x)$, for each $\la\in\mathbb{C}$, the unique solution  of problem \eqref{subpro} (defined below) $f(x;\la)$  and its derivative $f'(x;\la)$ are functions  of class $\mathcal{C}^{1}$ on $]0,L[\times \mathbb{C}$ (see Peano's theorem in \cite{hartman1964ordinary} 5.3).
Thus, for each continuity point of $\phi(x,\la)=p(x)f'(x;\la)/f(x;\la)$ we get that $\phi(x,\la)$ is of class $\mathcal{C}^{1}$.

Now, we define the following functions
\begin{equation}
\label{funciones}
h_{1}(\la)=\frac{f'(L;\la)}{\la f(L;\la)},\qquad h_{2}(\la)=\beta_{1}-\frac{\beta_{2}}{\la}.
\end{equation}
We remark that the regularity of the function is independent of the constant $C_{1}$.

Without loss of generality, we can consider a $C_{1}$ such that $f(0;\la)=1$.
Under that consideration, for $f(x;\la)$  the solution of the ODE in \eqref{problema} and the boundary condition at $x=0$, we get that is the unique regular solution for the following initial value problem
\begin{equation}
\label{subpro}
\left\lbrace \begin{array}{rcl}
                (p(x)f)'+(q(x)-\la r(x))f   &=&0\qquad 0<x<L\\
                                          f(0) &=&1\\
                                          f'(0)&=&\al_{1}\la +\al_{2}

             \end{array}
      \right.
\end{equation}
Reciprocally, the solution of \eqref{subpro} satisfies the ODE and the boundary condition at $x=0$ in our problem \eqref{problema}.

We note that function $h_{1}(\la)$ is well defined as long as it satisfies  $f(L;\la)\neq 0$.
Therefore, in order to obtain the domain of function $h_{1}(\la)$ we consider the following auxiliary boundary value problem.
\begin{equation}
\label{aux1}
\left\lbrace\begin{array}{rcl}
 \dis{(p(x)y')'-(q(x)- \lambda r(x) )y}&=&0\qquad 0<x<L \\
                                  y'(0)&=&(\alpha_{1}\lambda +\alpha_{2})y(0)\\
                                   y(L)&=&0
            \end{array}\right.
\end{equation}
Unlike in our problem, in  \eqref{aux1} there is a fix boundary condition, which eases (in part) the calculations, since some of the results on regular SLP are directly applicable, e.g. Sturm's comparison theorems.
For this problem, we have the following result on the spectrum.
\begin{te}
\label{secundario}Under the condition  $\alpha_{1}>0$, the  spectrum of \eqref{aux1} is an  ordered set
$$\eta_{-0}<0<\eta_{0}<\eta_{1}<\eta_{2}<\dots\nearrow \infty$$
with the property that  the eigenfunction $y(x;\eta_{l})$ has exactly $\vert l\vert$ zeros in $]0,L[$.
\end{te}
The proof is presented in \cite{binding2004transformation}. In Appendix \ref{demoaux} we present an alternative proof.

Using this result, we note that function$h_{1}(\la)$ is well defined if $\la\neq \eta_{n}$, with $\eta_{n}$ an eigenvalue of \eqref{aux1} and $\la\neq 0$.
Now, the analysis of the monotone behavior of function $h_1(\la)$ is developed on the branches
\begin{equation}
\label{ramas}
\mathcal{B}_{-1}=]-\infty,-\eta_{0}[;\quad\mathcal{B}_{-0}=]-\eta_{0},0[;\quad \mathcal{B}_{0}= ]0,\eta_{0}[;\quad \dots \quad \mathcal{B}_{n}= ]\eta_{n-1},\eta_{n}[
\end{equation}
This way, on each branch $\mathcal{B}_{n}$  we study the solutions of equation \eqref{laecuacion} through an analysis of the function graph intersection of $h_{1}(\la)$ and $h_{2}(\la)$.
Knowing that $\beta_{2}>0$ we get that $h_{2}(\la)$ is increasing on each branch $]-\infty,0[$ and $]0,\infty[$.
The monotony of $h_{1}(\la)$ is obtained through the analysis of an ODE in terms of $h_{1}(\la)$ and the coefficient functions $p(x)$, $q(x)$ and $r(x)$.
For this point, the regularity of functions $f(L;\la)$ and $f'(L;\la)$ respect to $\la$ is fundamental.
In the following subsection we cover the subject of regularity.

\begin{rem}
For the rest of this article we will denote by $f(x;\la)$ a regular function satisfying the following equations
\begin{equation}
\label{main2}
\left\lbrace\begin{array}{rcl}
 \dis{(p(x)f')'-(q(x)- \lambda r(x))f}&=&0\qquad 0<x<L \\
                                          f'(0)&=&(\alpha_{1} \lambda +\alpha_{2})f(0)
            \end{array}\right.
\end{equation}
Moreover, for all class of results on the regularity of ODEs, we consider  $f(x;\la)$ the solution of \eqref{main2} such that $f(0;\la)=1$.
As we have already noted, this consideration does not yield a loss of generality, since it doesn't modify the behavior of function $h_{1}(\la)$ in \eqref{funciones}.
\end{rem}

\subsection{On the existence of $\la_{0}$ and $\la_{-0}$}\label{valorespropios}
In this part we  prove the existence of eigenvalues $\la_{-0}<0<\la_{0}$ such that the respective eigenfunctions don't change sign in the interval $]0,L[$.
In the proof of this existence result, we present the main idea to prove the monotonicity of  $h_{1}(\la)$ in \eqref{funciones}.
\begin{lemma}
\label{primervalorproio}
There exist $\la_{-0}$  and $\la_{0}$ two eigenvalues of \eqref{problema} such that $\eta_{-0}<\la_{-0}<0<\la_{0}<\eta_{0}$ and the respectively eigenfunctions are positives in $]0,L[$.
\end{lemma}
In the proof of this Lemma we obtain an ODE in terms of  $h_{1}(\la)$.
Analyzing the sign of  $\frac{d}{d\la}h_{1}$ we get that $h_{1}(\la)$ is monotone decreasing for $\la\in]0,\eta_{0}[$.
The proof for the existence of $\la_{-0}$ is similar.

The argument presented in this part will be used for our proof of the existence of a sequence of eigenvalues in \eqref{orden}.
For the general case, other class of technical arguments are required.

\noindent{\bf Proof:} We begin considering
\begin{equation}\label{va}
\va(x,\la)=p(x)f'(x;\la)/f(x;\la),
\end{equation}
with $f(x;\la)$ solution of \eqref{main2}.
From Peano's theorem (see Hartman \cite{hartman1964ordinary} 5.3) we know that $\va(x,\la)$ is of class $\mathcal{C}^{1}$ in $]0,L[\times\cc$.

It follows directly that
\begin{equation*}
\begin{array}{rcl}
\displaystyle{\frac{\pp \va }{\pp x}}&=&\displaystyle{\frac{(pf')'}{f}-\mu_{0}\frac{(f')^{2}}{f^{2}}}\\
                      &=& \displaystyle{-\la r(x)+q(x)-\frac{1}{p(x)}\va^{2}}.
\end{array}
\end{equation*}
Now, taking the partial derivative with respect to $\lambda$ and interchanging the order of the derivatives, we get
\begin{equation*}
\dis{\frac{\pp }{\pp x}\frac{\pp \va }{\pp \la}=-r(x)-2\frac{\va}{p(x)}\frac{\pp \va }{\pp \lambda}}.
\end{equation*}
Taking $\dis{\psi(x,\la)=\frac{\pp}{\pp\la}\va}$, we rewrite the above equation as follows
$$ \dis{ \frac{\pp \psi}{\pp x} +2\frac{f'(x;\la)}{f(x;\la)}\psi =-r(x)      }  $$
and, therefore we get
\begin{equation}
\frac{\pp}{\pp x}\left(f^{2}\psi\right)=-rf^{2}.
\end{equation}
Thus, integrating between $[0,L]$ we obtain the following relation
\begin{equation}
\label{rel1}
f^{2}(L;\la)\psi(L;\la)=-\int_{0}^{L} f^{2}(x;\la) r(x)dx+f^{2}(0;\la)\psi(0;\la)
\end{equation}
On the other hand, from the ODE in \eqref{main2} we know that
\begin{equation}
\label{integral}
\begin{array}{rcl}
\displaystyle{-\int_{0}^{L} r(s)f^{2}(s;\la)ds }&=&\dis{\frac{1}{\la}\left\lbrace \int_{0}^{L} f(pf')'ds-\int_{0}^{L} q f^{2}ds  \right\rbrace  }\\
 &=&\dis{\frac{1}{\la}\left\lbrace f(L;\la)p(L)f'(L;\la)-f(0;\la)p(0)f'(0;\la)-\right.}\\
 & & \displaystyle{\left.-\int_{0}^{L} p(f')^{2}ds-\int_{0}^{L} q f^{2}ds  \right\rbrace  }\\
 &=&\dis{\frac{1}{\la}\left\lbrace f(L;\la)^{2}\va(L,\la)-f^{2}(0;\la)(\al_{1}\la+\al_{2})-\right.}\\
 & &\dis{-\left.\int_{0}^{L}\left[ p(f')^{2}+ q f^{2}\right]ds  \right\rbrace  }\\
\end{array}
\end{equation}
Using the expression for the boundary condition at $x=0$ it follows that $\psi(0,\la)=\al_{1}$ and, therefore  we get
$$ f^{2}(L;\la)\psi(L;\la)-\frac{1}{\la} f^{2}(L;\la) \va(L,\la)= -f^{2}(L;\la)\frac{\al_{2}}{\la}-\dis{\frac{1}{\la}\left\lbrace \int_{0}^{L}\left[ p(f')^{2}+ q f^{2}\right]ds  \right\rbrace  }  $$
Finally, knowing that $\psi =\frac{\pp \va}{\pp \la}$, from the above equation we obtain
\begin{equation}
\label{ODE}
\frac{d}{d\la}\left(\frac{1}{\la}\va(L,\la)\right)=-\frac{1}{(\la f(L;\la))^{2}}\left\lbrace \al_{2}+\int_{0}^{L}\left[ p(f')^{2}+ q f^{2}\right]ds  \right\rbrace
\end{equation}
Hence, the function $h_{1}(\la)=\frac{1}{\la}p(L)\frac{f'(L;\la)}{f(L;\la)}$ is a decreasing function for $\la>0$ while $f(L;\la)>0$.
The Next step is to study the equation \eqref{laecuacion} through the graph  intersection of $h_{1}(\la)$ and $h_{2}(\la)=\beta_{1}(k)-\frac{\beta_{2}(k)}{\la}$.

Knowing that $\al_{1}$ and $\al_{2}$ are positive, from the consideration $f(0;\la)=1$ we get that $f'(0;\la)>0$ for each $\la\in\mathcal{B}_{0}=]0,\eta_{0}[$.
On the other hand, as $\la>0$, from the maximum principle for second order linear differential ODE, we get $f'(L;\la)<0$.
Therefore, from the definition of $h_{1}(\la)$ and $\eta_{0}$, the principal eigenvalue of the auxiliary non-regular SLP \eqref{aux1}, we obtain
$$\lim_{\la\searrow 0}h_{1}(\la)=\infty;\qquad \lim_{\la\nearrow\mu_{0}}h_{1}(\la)=-\infty$$
Thus, $h_{1}:\mathcal{B}_{0}\to\rr$ is a surjective monotonic decreasing function.
Finally, as $\beta_{2}>0$ we get that $h_{2}(\la)$ in \eqref{funciones} is increasing function in $\mathcal{B}_{0}$ and therefore there exists one graph intersection between $h_{1}(\la)$ and $h_{2}(\la)$ in $\mathcal{B}_{0}$.
If we denote by $\la_{0}$ the unique solution of \eqref{laecuacion} in $\mathcal{B}_{0}$, from the definition of $h_{2}(\la)$ and knowing that $f(x;\la_{0})$ satisfy \eqref{main2}, we get that $\la_{0}$ is an eigenvalue of \eqref{problema} where $f(x;\la_{0})>0$ in $]0,L[$.
Using similar argument we obtain the existence of $\la_{-0}\in]\eta_{-0},0[$, eigenvalue of \eqref{problema}, such that $f(x;\la_{-0})>0$ in $]0,L[$.
$\square$
\begin{rem}
We use a similar argument for the description of the other eigenvalues in \eqref{orden}.
For the general case we need other technical considerations, but our proof is centered in the sign of $\frac{d}{d\la}h_{1}(\la)$.
For this end, in the next section we present some results on the behavior of the zeros of $f(x;\la)$ (solution of \eqref{main2}).
\end{rem}

\subsection{Regular SLP associated with \eqref{problema}}

In this subsection we obtain a regular SLP where its spectrum is related with the spectrum of our problem \eqref{problema} by means of a Crum-Darboux type transformation. We use the ideas and terminology given in \cite{binding2004transformation}.
From the Lemma \ref{primervalorproio}, let $\la_{0}>0$ be an eigenvalue of \eqref{problema} such that $y_{0}(x)=f(x;\la_{0})>0$.

Consider the following Crum-Darboux type of transformation
\begin{equation}
\label{transformacion}
g(x)=p(x)f'(x;\la)-p(x)f(x;\la)\frac{y_{0}'(x)}{y_{0}(x)}
\end{equation}
with $y_{0}(x)=f(x;\la_{0})$, where $f(x;\la)$ is solution of \eqref{main2}.\\
Therefore,  $g(x)$ in (\ref{transformacion}) satisfy the following regular boundary value problem
\begin{equation}
\label{problem2}
(P')\left\lbrace\begin{array}{rcl}
 \dis{(\tilde{p}(x)g')'-(\tilde{q}(x)-\la\tilde{r}(x))g}                &=&0\qquad0<x<L \\
                                               g'(0)&=&-\tilde{\al}  g(0)\\
                                               g'(L) &=&-\tilde{\beta}g(L)
            \end{array}\right.
\end{equation}
where
$$\tilde{p}=\frac{1}{r(x)};\quad \tilde{r}(x)=\frac{1}{p(x)}\la_{0}-\left( -\frac{1}{r}\frac{y_{0}'}{y_{0}} \right)'+\frac{1}{r}\left(\frac{y'_{0}}{y_{0}}\right)^2;\quad  \tilde{r}(x)=\frac{1}{p(x)}$$
and
\begin{equation}
\label{condicionesregular}
\tilde{\al}_{1}=\frac{r(0)}{\al_{1}} +\frac{y_{0}'(0)}{y_{0}(0)};\qquad \tilde{\beta}_{1}=\frac{r(L)}{\beta_{1}} +\frac{y_{0}'(L)}{y_{0}(L)}
\end{equation}

The deduction of \eqref{problem2} is presented in Appendix  \ref{regularslp}.

We can derive directly the following result on the spectrum of \eqref{problema}:
\begin{coro}
\label{soloreales}
The spectrum of \eqref{problema} is a countable real set.
\end{coro}
\noindent{\bf Proof:} The proof is obtained directly, by contradiction, considering the regular SLP \eqref{problem2}.
$\square$

\subsection{Monotonic behaviors of the zeros of $f(x;\la)$.}

In this part, we study the behavior of the zeros of $f(x;\la)$, a solution of the initial value problem \eqref{main2}.
Through a collection of lemmas, we present some fundamental tools for our analysis on the behavior of the function $h_{1}(\la)$ defined in \eqref{funciones}.
We note that in the case where $\al_{1}<0$ and $\beta_{1}>0$, using Sturm's comparison theorems we have a decreasing behavior of $h_{1}(\la)$ for $\la>0$.
In the case where $\al_{1}\beta_{1}\geq 0$, this results can be adapted using Majorant and Minorant Sturm's problems (see Hartman \cite{hartman1964ordinary} 11.3). In Appendix \ref{otroespectro} we present the analysis for the case $\al_{1}\beta_{1}\geq 0$.

In the following, we consider  $z_{j}(\la)$  to denote the $j-$th zero of $f(x;\la)$,  i.e., $f(x;\la)$ has exactly $j-1$ zeros in $]0,z_{j}(\la)[$, also satisfying $f(z_{j}(\la);\la)=0$.
We describe  $z_{j}(\la)$ as a regular function of the variable $\la$ by means of the implicit function.

For the development of our argument, we begin by noting the regularity of $f(x;\la)$, respective to the parameter $\la$.

Since $r(x)$ and $q(x)$ are positive, through oscillatory results for second order ODE, we know that there exist a $\la^{*}>0$ such that  $ f(x;\la^{*})=0 $ has at least one root in $]0,L[$.

Considering $\la^{*}$ fixed, we denote by $x^{*}$ the first zero of $f(x;\la^{*})$ in  $]0,L[$, i.e., we get that $f^{2}(x;\la^{*})>0$ in $]0,x^{*}[$ and $f(x^{*};\la^{*})=0$.
On the other hand, if we assume that $f'(x^{*};\la^{*})=0$, considering the initial value problem
\begin{equation}
\label{cauchy}
(p(x)f')'-(q(x)-\la^{*} r(x))f=0,\qquad f(x^{*})=f'(x^{*})=0
\end{equation}
we have that $f(x;\la^{*})=0$, obtaining a contradiction.
Thus $f'(x^{*};\la^{*})\neq 0$.

Now, consider $\phi:[0,L]\times \rr\to\rr$ defined by $\phi(x,\la)=f(x;\la)$.
From the regularity of $p(x),\ q(x)$ and $r(x)$, using the Peano's theorem (see Hartman \cite{hartman1964ordinary} 5.3), for each neighborhood $\Omega^{*}$ of $(x^{*},\la^{*})$ we have that $\phi\in\mathcal{C}^{1}(\Omega)$.
Moreover, as $f'(x^{*};\la^{*})\neq 0$, from the definition of $\phi(x,\la)$ we also have that $\dis{ \frac{\pp }{\pp x}\phi(x^{*},\la^{*})\neq 0 }$.
Therefore, from the implicit function theorem (IFT), we know that there exist $I_{1}$ and $I_{2}$, neighborhood of $\la^{*}$  and $x^{*}$ respectively, and a unique $\mathcal{C}^{1}$ function $ g:I_{1}\to I_{2}$, such that
$$ \phi(g(\la);\la)=0,\qquad \textrm{for each}\ \la\in I_{1}.$$
Moreover, from the IFT it follows that
$$  \frac{dg}{d\la}=-\dis{\frac{\frac{\pp }{\pp \la}\phi(g(\la),\la)  }{\frac{\pp }{\pp x}\phi(g(\la),\la)}   }   $$
Thus, as $\phi(x,\la)=f(x;\la)$, given some $\la\in\cc$, such that $f(x;\la)=0$ has at least one root in $]0,L[$, there exists some neighborhood $I_{1}$ of $\la$ such that the first root of $f(x;\la)=0$ can be defined as a regular function $z_{1}:I_{1}\to \rr$ satisfying
\begin{equation}
\label{signo2}
\dis{ \frac{d z_{1}}{d\la}= -\dis{\frac{\frac{\pp }{\pp \la}f(z_{1}(\la);\la)  }{f'(z_{1}(\la);\la)}   }    }.
\end{equation}
This argument can be repeated for each zero denoted by $z_{j}(\la)$.

The aim is to obtain the sign of the derivatives of $z_{j}(\la)$.
We note that, through the oscillatory results, it's possible to prove that if $l\geq j$, and $\la \in\mathcal{B}_{l}$, then $z_{j}(\la)$ is well defined and its domain is given by $]\eta_{j},\infty[$ with $\eta_{j}$ eigenvalue of \eqref{aux1}.

We begin proving that for each $\la<\eta_{-0}$ the function $f(x;\la)$ has exactly one zero in $]0,L[$.
\begin{lemma}
\label{numerodezeronegativos}
For each $\la\in]-\infty,\eta_{-0}[$ the solution of (\ref{main2}) has exactly one zero in $]0,L[$.
\end{lemma}
\noindent{\bf Proof:} For $\la<\eta_{-0}$ we have that $q(x)-\la r(x)>0$ and using oscillatory results, $f(x;\la)$ is non oscillatory in $]0,L[$.
Thus, $f(x;\la)$ has at most one zero in $]0,L[$.

Now, we use Sturm's first comparison theorem using a Sturm majorant at (\ref{main2}).
For $a=\min\lbrace q(x): 0\leq x\leq L \rbrace$; $b=\min\lbrace r(x): 0\leq x\leq L \rbrace$ and $c=\min\lbrace p(x): 0\leq x\leq L \rbrace$, consider
\begin{equation}
\label{majo2}
\left\lbrace\begin{array}{rcl}
 \dis{(c\tilde{f}')'-(a- \lambda b)\tilde{f}}&=&0\qquad 0<x<L \\
                                          \tilde{f}'(0)&=&(\alpha_{1} \lambda +\alpha_{2})\tilde{f}(0)
            \end{array}\right.
\end{equation}
Directly, from the definition of $a,\ b$ and $c$, we have that \eqref{majo2} is a majorant problem at \eqref{main2}.
Under our consideration that $p'(x)>0$, we have that $c=p(0)$.

On the other hand, taking $\overline{c}=\sqrt{(a-b\lambda)/c}$, we get that
\begin{equation}
\label{somajorant}
\tilde{f}(x;\la)=\cosh(\overline{c}(\la)x)+\frac{(\al_{1}\la+\al_{2})}{\overline{c}(\la)}\sinh(\overline{c}(\la)x)
\end{equation}
defines a solution of  (\ref{majo2}).
As $\lim_{\la\to-\infty}\tanh(\overline{c}(\la)(s+L))=1$ and  $\overline{c}(\la)\sim O(1/\sqrt{-\la})$ when $\la\to-\infty$ then, there exists some $\la^{*}$ such that the equation $\tilde{f}(x;\la^{*})=0$ has one root in $]0,L[$.

As (\ref{main2}) and (\ref{majo2}) have the same boundary condition in $x=0$, we have
\begin{equation}
\label{secondsturm1}
\frac{c \tilde{f}'(0;\la)}{\tilde{f}(0;\la)}\geq \frac{p(0)f'(0;\la)}{f(0;\la)}.
\end{equation}
Finally, using the Sturm's first comparison theorem, we have that $f(x;\la)$ has exactly one zero in $]0,L[$.
   $\square$

\begin{rem}
Through the analysis of the behaviors of $z_{j}(\la)$, we develop our arguments to conclude that $h_{1}(\la)$ is a decreasing function on each of its branches $\mathcal{B}_{n}$ defined in \eqref{ramas}.

We remark that, when $\la\in]\eta_{-0},\eta_{0}[$ we have $f(x;\la)>0$ in $]0,L[$.
\end{rem}

The analysis in the cases $\la<\eta_{-0}$ and $\la> \eta_{1}$ will be studied separately in the following two Lemmas.
We begin with $\la<\eta_{-0}$ as follows:
\begin{lemma}
\label{zeros2}
 For $\underline{\la}<\overline{\la}<\eta_{-0}$, we get $z_{1}(\underline{\la})<z_{1}(\overline{\la})$.
\end{lemma}
\noindent{\bf Proof:} The existence of $z_{1}(\la)$ for $\la<\eta_{-0}$ was proved in Lemma \ref{numerodezeronegativos}.
Let $\dis{ \va(x,\la)=p(x)\frac{f'(x;\la)}{f(x;\la)} }$ be a regular function in $]0,z_{1}(\la)[$.

Now, taking
$$\psi(x;\la)=\dis{\frac{\pp \va}{\pp\la}}=p(x)\frac{ \dis{\frac{\pp f'}{\pp \la} f - f' \frac{\pp f}{\pp \la} } }{ f^{2} (x;\la)   },$$
we get
\begin{equation}
\label{limite}
\dis{\lim_{x\nearrow z_{1}(\la)}f^{2}(x;\la)\Psi(x;\la)=-p(z_{1}(\la))f'(z_{1}(\la);\la)\frac{\pp}{\pp\la}f(x;\la)\left\vert_{x=z_{1}(\la)}\right.
  }
\end{equation}
Thus, in the following steps we commit the limit procedure.
From the ODE in \eqref{problema} we know that
\begin{equation}
\label{eq reducida 2}
\frac{\pp}{\pp x}\left( f^{2}(x;\la) \psi(x,\la)   \right)=-r(x)f^{2}(x;\la).
\end{equation}
On the other hand, from $f(x;\la)=1$ and $f'(0;\la)=\al_{1}\la+\al_{2}$, we get that $\va(0,\la)=p(0)(\al_{1}\la +\al_{2})$, and therefore
$$\psi(0,\la)=\dis{\frac{\pp\va(0;\la)}{\pp\la} =p(0)\frac{\pp}{\pp\la}(\al_{1}\la+\al_{2}) =p(0)\al_{1}}$$
Thus, using similar argument as in the limit presented in (\ref{limite}), integrating the left hand side of the equation (\ref{eq reducida 2}) between $0$ and $z_{1}(\la)$, we obtain
\begin{equation}
\label{izq}
\begin{array}{rcl}
\dis{\int_{0}^{z_{1}(\la)}\frac{\pp}{\pp x}\left( f^{2}(x;\la) \psi(x,\la)   \right)dx }&=& \dis{-p(z_{1}(\la))f'(z_{1}(\la);\la) \frac{\pp}{\pp\la} f(z_{1}(\la);\la)-f^{2}(0;\la)\psi(0,\la)}\\
                             &=&\dis{-p(z_{1}(\la))f'(z_{1}(\la);\la) \frac{\pp}{\pp\la} f(z_{1}(\la);\la)-p(0)\al_{1} }.
\end{array}
\end{equation}
Now, integrating the right side of (\ref{eq reducida 2}) between $0$ and $z_{1}(\la)$ we obtain
\begin{equation}
\label{der}
\begin{array}{rcl}
\dis{-\int_{0}^{z_{1}(\la)}r(x)(f(x;\la))^{2}dx}&=& \dis{ \frac{1}{\la} \int_{0}^{z_{1}(\la)} f(pf')'dx -  \frac{1}{\la} \int_{0}^{z_{1}(\la)} qf^{2}dx}\\
                                               &=& \dis{\frac{1}{\la} \left((fpf')(z_{1}(\la))-(fpf')(0) - \int_{0}^{z_{1}(\la)} (p(f')^{2}+qf^{2})dx \right)}\\
                                               &=&\dis{-\frac{p(0)}{\la}(\al_{1}\la +\al_{2})  -\frac{1}{\la}\int_{0}^{z_{1}(\la)} (p(f')^{2}+qf^{2})dx  }

\end{array}
\end{equation}
Thus, if we integrate (\ref{eq reducida 2}) between $0$ and $z_{1}(\la)$, from (\ref{izq}) and (\ref{der}) we obtain
$$ \dis{-p(z_{1}(\la))f'(z_{1}(\la);\la) \frac{\pp}{\pp\la} f(z_{1}(\la);\la)-p(0)\al_{1}}=\dis{-\frac{p(0)}{\la}(\al_{1}\la +\al_{2})  -\frac{1}{\la}\int_{0}^{z_{1}(\la)} (p(f')^{2}+qf^{2})dx  }  $$
and therefore
\begin{equation}
\label{casi}
\dis{p(z_{1}(\la))f'(z_{1}(\la);\la) \frac{\pp}{\pp\la} f(z_{1}(\la);\la )}=\dis{\frac{1}{\la}\left(\al_{2} + \int_{0}^{z_{1}(\la)} (p(f')^{2}+qf^{2})dx    \right)  }
\end{equation}
On the other hand, from the IFT we get \eqref{signo2}, i.e.$\dis{\frac{d z_{1}(\la)}{d\la} =-\frac{ \frac{\pp}{\pp\la} f(z_{1}(\la);\la )}{f'(z_{1}(\la);\la)} }$ and therefore the sign of $\frac{dz_{1}}{d\la}$ can be obtained from (\ref{casi}).
Finally, as $\la<\eta_{-0}<0$ we get $\dis{\frac{d z_{1}(\la)}{d\la}>0 }$.
  $\square$

Now, we analyze the behavior of $z_{n}(\la)$ for $\la>\eta_{0}$.
\begin{lemma}
\label{zeros1}
Let $\la>\eta_{0}$ and $z_{n}(\la)$ the $n-$th zero of $f(x;\la)$.
Thus, $\frac{d z_{n}}{d\la}<0.$
\end{lemma}
\noindent{\bf Proof:} Assume that  $f(x;\la)$ has at least $n$ zeros in $]0,L[$, denoted by $z_{1}(\la)<z_{2}(\la)<\dots<z_{n}(\la)$.
For $z_{1}(\la)$, using \eqref{casi} we get that $\frac{d z_{1}}{d\la}<0$.
Now, integrating (\ref{eq reducida 2}) between $z_{1}(\la)$ and $z_{2}(\la)$, we get
\begin{equation}
\label{paso1}
-p(z_{2})f'(z_{2};\la)\frac{\pp}{\pp\la}f(z_{2};\la)+p(z_{1})f'(z_{1};\la)\frac{\pp}{\pp\la}f(z_{1};\la)=-\dis{\int_{z_{1}}^{z_{2}} r f^{2}dx}.
\end{equation}
Using the implicit function theorem, from \eqref{paso1} we obtain
\begin{equation}
\label{paso2}
p(z_{2})(f'(z_{2};\la))^{2}\frac{dz_{2}}{d\la}=p(z_{1})(f'(z_{1};\la))^{2}\frac{dz_{1}}{d\la} -\dis{\int_{z_{1}}^{z_{2}} r f^{2}dx}
\end{equation}
Knowing that $\frac{dz_{1}}{d\la}<0$, from \eqref{paso2} we get $\frac{dz_{2}}{d\la}<0$.
Following similar steps we obtain the proof for any $z_{j}(\la)$ with $j=3,\dots,n$.
$\square$

\subsection{Behavior of $h_{1}(\la)$ as $\la\in]-\infty,\eta_{-0}[$.
Asymptotic representation.}\label{asym}
Using the asymptotic representation for  $h_{1}(\la)$ presented in this section, we obtain results of existence of $\la_{-1}<0$, eigenvalue of \eqref{problema} such that the corresponding eigenfunctions change sign exactly once in $]0,L[$ (see Lemma \ref{numerodezeronegativos}).

We continue working with $f(x;\la)$, solution of \eqref{main2}.
In the next result we consider $h_{1}(\la)$ defined in \eqref{funciones} and $\eta_{-0}<0$ eigenvalue of \eqref{aux1} such that the respective eigenfunction satisfy $y(x;\eta_{0})>0$ in $]0,L[$.
\begin{lemma}
\label{sobrelafuncion}
The function $h_{1}:]-\infty,\eta_{-0}[\to]-\infty,0[$ is a bijective map and monotonic decreasing function and its asymptotic representation is given by $h_{1}(\la)=O(\la^{-1/2})$ as $\la\to-\infty$.
\end{lemma}
\begin{rem}
Directly from the fact stated in the previous lemma, we note that in case the $\beta_{1}<0$, the graphs of functions $h_{1}(\la)$ and $h_{2}(\la)$ do not intersect in $]-\infty,\eta_{-0}[$, therefore $\la_{-1}$ doesn't exist in this case.
\end{rem}

\noindent{\bf Proof:} We begin by obtaining the sign of $h_{1}(\la)$.
Since $q(x)-\la r(x)>0$, when $\la\leq 0$, knowing that  $\la<\eta_{-0}<-\frac{\al_{2}}{\al_{1}}$, from  maximum principle for second order linear ODE (see \cite{protter1984maximum} Chapter 1) we have that $f'(L;\la)\cdot f(L;\la)>0$.
Thus, $h(\la)<0$ for $\la<\eta_{-0}$.
Moreover, we have $\lim_{\la\nearrow \eta_{-0}}h_{1}(\la)=-\infty$.

For the asymptotic representation as $\la\to\infty$, we use the notion of Sturm  majorant and Sturm  minorant problem for \eqref{main2}.

We begin considering the Sturm majorant at \eqref{main2} given in \eqref{majo2}.
Knowing that $\tilde{f}(x;\la)$ in \eqref{somajorant} has one zero for $\la<\eta_{-0}$, from the condition at $x=0$ and Lemma \eqref{numerodezeronegativos}, through Sturm's second comparison theorem, we get
\begin{equation}
\label{upperbound}
c\frac{\tilde{f}'(L;\la)}{ \tilde{f}(L;\la)}>p(L)\frac{f'(L;\la)}{f(L;\la)}
\end{equation}
To obtain a lower bound for the asymptotic behavior of $h_{1}(\la)$, we use \eqref{upperbound} considering that $\tilde{f}'(L;\la)/\tilde{f}(L;\la)\sim \overline{c}(\la)\coth(\overline{c}L)$, as $\la\to-\infty$.\\
Thus, knowing that $\coth(\overline{c}(\la)L)\sim 1$ as $\la\to -\infty$, from \eqref{upperbound} we obtain
\begin{equation}
\label{cota1}
-\sqrt{b\cdot c}(-\la)^{-1/2}<h_{1}(\la),\qquad\textrm{as}\ \la\to-\infty
\end{equation}
Now, we use a Sturm minorant problem at \eqref{main2} and  following a  similar argument we obtain an upper bound.
Consider  $\underline{c}(\la)=\sqrt{ (A-B\la)/C }$, with  $A=\max\lbrace q(x): 0\leq x\leq L \rbrace$; $B=\max\lbrace r(x): 0\leq x\leq L \rbrace$ and $C=\max\lbrace p(x): 0\leq x\leq L \rbrace$.
Using this upper bounds we define the following Sturm  minorant problem:
\begin{equation*}
\label{mino}
\left\lbrace\begin{array}{rcl}
 \dis{y''-(\underline{c}(\la))^{2}y}&=&0\qquad0<x<L \\
                                          y(L)&=&0
 \end{array}\right.
\end{equation*}
Following a similar argument as in the proof of the bound given in (\ref{cota1}) we obtain
\begin{equation}
\label{cota2}
h_{1}(\la)>C\frac{\underline{c}(\la)}{\la}\textrm{cotanh}(-\underline{c}(\la)L).
\end{equation}
Hence, since $\underline{c}(\la)$ and $\overline{c}(\la) $ have the same asymptotic behavior as $\la\to-\infty$, given by $\sqrt{-\la}$, using the bounds in (\ref{cota1}) and (\ref{cota2}), we obtain that $\tilde{h}(\la)=O((-\la)^{-1/2})$ as $\la\to-\infty$.
Finally, as $h_{1}(\la)$ is a decreasing function, we get that $h_{1}:]-\infty,\eta_{-0}[\to]-\infty,0[$ is a surjective function.
$\square$

\section{Proofs of  Theorems \ref{principal} and \ref{principal2}.}\label{demoprincipales1}
In this section we develop the proofs of Theorems \ref{principal} and \ref{principal2} through the study of function $h_{1}(\la)$.
The collection of eigenvalues in \eqref{orden} is obtained though the analysis of graph intersection of the functions $h_{1}(\la)$ and $h_{2}(\la)$.
The monotone behavior of function $h_{1}(\la)$ is fundamental for this part.
The main idea of the proof was presented inf the proof of Lemma \ref{primervalorproio}.
We use Lemmas \ref{zeros2} and \ref{zeros1} for the oscillatory results of the eigenfunctions in Theorem  \ref{principal} and for the separation Theorem \ref{principal2}.
We begin by describing the eigenvalues of  \eqref{orden} and continue working with $f(x;\la)$, solution of \eqref{main2}.

\noindent{\bf Proof Theorem \ref{principal}} From Corollary \ref{soloreales} we know that the spectrum must be a real subset.
Then, through graph intersection analysis of functions $h_{1}(\la)$ and $h_{2}(\la)$ we obtain the result.

Consider $\la\in\mathcal{B}_{n}=]\eta_{n-1},\eta_{n}[$, with $n=1,2,\dots$.
From the boundary condition at $x=0$ of \eqref{aux1}, through the Sturm's first comparison theorem we get $f(x;\la)$ has at least $l$ zeros in $]0,L[$.
Let $z_{l}(\la)$ the $l-$th zero of $f(x;\la)$ as in Lemma \ref{zeros1}.
Consider that $l\geq n$ and assume that $f(x;\la)$ doesn't change sign in $]z_{l}(\la),L[$.\\
Similar to the proofs of Lemmas \ref{primervalorproio} and \ref{zeros2}, we have
\begin{equation}
\label{1}
f^{2}(L;\la)\frac{\pp}{\pp\la}\va(L,\la)=-p(z_{l}(\la))f'(z_{l}(\la);\la)\frac{\pp}{\pp\la}f(z_{l}(\la);\la)-\int_{z_{l}(\la)}^{L}k^{2}p'f^{2}d\tau,
\end{equation}
with $\va(L,\la)=p(L)f'(L;\la)/f(L;\la)$.
As $f(z_{l}(\la);\la)=0$, using similar limit argument as in \eqref{limite}, we obtain
\begin{equation*}
\begin{array}{rcr}
f^{2}(L;\la)\left( \frac{d}{d\la}\va(L,\la) -\frac{1}{\la}\va(L,\la) \right)&=&-p(z_{l}(\la))f'(z_{l}(\la);\la)\frac{\pp}{\pp\la}f(z_{l}(\la);\la) -\\
\\
& &\dis{-\frac{1}{\la}\int_{z_{l}(\la)}^{L}p(f'^{2}+k^{2}f^{2})d\tau}
\end{array}
\end{equation*}
From the definition of $\va(L;\la)$, the above identity can be written as follow
\begin{equation}
\label{sigder}
\begin{array}{rcr}
\dis{\frac{d}{d\la}\left( \frac{\va(L,\la)}{\la}   \right)}&=&\dis{-\frac{1}{f^{2}(0;\la)}\left\lbrace \frac{f'(z_{l}(\la);\la)\frac{\pp}{\pp\la}f(z_{l}(\la);\la)}{\la}+\right.}\\
               & & \dis{\left.+\frac{1}{\la^{2}}\int_{z_{l}(\la)}^{0}p(f'^{2}+k^{2}f^{2})d\tau\right\rbrace}.
\end{array}
\end{equation}
On the other hands, from the IFT we have that
$$ \frac{d z_{l}(\la)}{d\la}=-\dis{ \frac{\frac{\pp}{\pp\la}f(z_{l}(\la);\la)}{f'(z_{l}(\la);\la)} } $$
Using Lemma \ref{zeros1} we know $\frac{d z_{l}}{d\la}<0$, from \eqref{sigder} we obtain that $h_{1}(\la)$ in (\ref{funciones}) is decreasing on each branch $\mathcal{B}_{n}=]\eta_{n-1},\eta_{n}[$, with $n=1,2,\dots$.
Moreover, as $\lim_{\la\searrow\eta_{n-1}}h_{1}(\la)=\infty$ and $\lim_{\la\nearrow\eta_{n}}h_{1}(\la)=-\infty$, we obtain that $h_{1}:\mathcal{B}_{n}\to\rr$ is a surjective function.
As $h_{2}(\la)$ is an increasing function, in each $\mathcal{B}_{n}$ there exist one graph intersection denoted by $\la_{n}$.
Knowing that $f(x;\la)$ is a solution of \eqref{main2}, from the definition of $h_{2}(\la)$ we get that $\la_{n}$ is a real eigenvalue of \eqref{problema}.\\
Now, we prove the oscillatory results for the eigenfunctions by reduction to absurd.
Assume that $l>n$, thus there exist $z_{n}(\la)<z_{l}(\la)$ zeros of $f(x;\la)$.
Taking
$$ \displaystyle{L=\lim_{\la\nearrow \eta_{n}}z_{n}(\la)< \lim_{\la\nearrow \eta_{n}}z_{l}(\la)< L,}$$
reaching a contradiction.
Therefore, the eigenfunction $f(x;\la_{n})$ Changes of sign exactly $n$ times in $]0,L[$.
Hence, we have the following collection of eigenvalues $\la_{1}<\la_{2}<\la_{3}<\dots$.
The existence of $\eta_{-0}<\la_{-0}<0<\la_{0}<\eta_{0}$ is presented in \eqref{primervalorproio}.
To finish our proof, we show the existence of $\la_{-1}<\eta_{-0}$.
From Lemma \ref{sobrelafuncion} we know that $\lim_{\la\to-\infty}h_{1}(\la)=0$.
Knowing that  $\lim_{\la\to-\infty}h_{2}(\la)=\beta_{1}<0$, through the monotony of functions $h_{1}(\la)$ and $h_{2}(\la)$ we get the existence and uniqueness of $\la_{-1}$, solution of \eqref{laecuacion}.
For the oscillatory result of $f(x;\la_{-1})$ we use Lemma \eqref{numerodezeronegativos}.

Finally, from the existence of $\la_{0}>0$ such that $f(x;\la_{0})>0$ in $]0,L[$, we have the result in Corollary \ref{soloreales}. Hence, all eigenvalues of \eqref{problema} are obtained as the graph intersections of $h_{1}(\la)$ and $h_{2}(\la)$ with $\la\in\rr$.
  $\square$

\noindent{\bf Proof Theorem \ref{principal2}:} Consider $\mathcal{B}_{l}$ in  (\ref{ramas}) with $l\geq 1$.
Given $\la\in\mathcal{B}_{l}$, from Theorem \ref{principal} we get that $f(x;\la)$ has exactly $l$ zeros in ]0,L[.
Now, consider $\la_{*}<\la^{*}$, both in $\mathcal{B}_{l}$.
From Lemma \ref{zeros1} we know that $z_{l}(\la^{*})<z_{l}(\la_{*})$.
Since the ODE in (\ref{pabstract}) is linear, we can assume that $\overline{f}(x)=f(x;\la^{*})$ and $\underline{f}(x)=f(x;\la_{*})$ are positive functions in $]z_{l}(\la^{*}),L]$ and $]z_{l}(\la_{*}),L]$ respectively.
For notational simplicity we consider $z_{l-1}=z_{l-1}(\underline{\la})$ and $z_{l}=z_{l}(\underline{\la})$.
Using Green's formula in $[z_{l-1},z_{l}]$, we obtain the following identity:
\begin{equation}
\label{signo}
p(z_{l})\overline{f}(z_{l})\underline{f}'(z_{l})-p(z_{l-1})\overline{f}(z_{l-1})\underline{f}'(z_{l-1})
=\displaystyle{(\la^{*}-\la_{*})\int_{z_{l-1}}^{z_{l}} r(\tau)\underline{f}(\tau)\overline{f}(\tau)d\tau}
\end{equation}
Knowing that $\underline{f}(x)>0$ in $]z_{l},L[$ and $\underline{f}(z_{l})=0$, then it holds that  $\underline{f}'(z_{j})>0$.
Similarly, knowing that $z_{j-1}<z_{j}$ are two consecutive zeros of $\underline{f}$, then $\underline{f}<0$ in $]z_{j-1},z_{j}[$ and $\underline{f}'(z_{j-1})<0$.

On the other hand, we know that $z_{j}(\la^{*})<z_{j}$, therefore it holds that $\overline{f}(z_{j})>0$.
Assuming that $z_{j}(\la^{*})<z_{j-1}$, we have that $\overline{f}(x)>0$ in $]z_{j-1},z_{j}[$, therefore we have $(\la^{*}-\la_{*})\int_{z_{l-1}}^{z_{l}} r(\tau)\underline{f}(\tau)\overline{f}(\tau)d\tau<0  $, reaching a contradiction with  \eqref{signo}.
Therefore, it holds that  $z_{j-1}<z_{j}(\overline{\la})<z_{j}$.
Following the same ideas, the proof for general case can be reached.
$\square$

\noindent{\bf Proof Corollary \ref{separacion}:} Using Theorem \ref{principal} we have that $f_{n}(x)$ has exactly $n$ zeros in $]0,L[$, which we denote by $z_{1}^{n}<z_{2}^{n}<\dots<z_{n}^{n}$.
Similar for $f_{n+1}$, considering $z_{1}^{n+1}<z_{2}^{n+1}<\dots<z_{n}^{n+1}<\dots<z_{n}^{n+1}$.

Using Lemma \ref{zeros1} we have that  $z_{1}^{n+1}<z_{1}^{n}$.
Now, using Theorem \ref{principal2} it holds that  $z_{1}^{n}<z_{2}^{n+1}<z_{2}^{n}$.
Hence $z_{1}^{n+1}<z_{1}^{n}<z_{2}^{n+1}$.
Proceeding similarly, the proof the general case can be reached.
$\square$


\section{Proof of Theorem \ref{asintotico}.
Asymptotic result for the spectrum}\label{demoprincipales2}
Now, we apply the Liouville transformation to obtain a normalized regular SLP.
Thus, the asymptotic behavior can be obtained through classical results for second order linear ODEs.

For
$$ P(x)=(\tilde{r}\tilde{p})^{1/2};\qquad G(x)=(\tilde{p}\tilde{r})^{-1/4}  $$
consider the following change of variable and transformation
\begin{equation}
\label{liouville}
t=\int_{0}^{x}P(u)du;\qquad g(x)=G(x)z(t)\qquad\textrm{(Liouville transformation)}.
\end{equation}
For
$$ R(t)=\left[\tilde{p}^{1/4}\tilde{r}^{-3/4}\frac{d}{dx}\tilde{p}\right]\frac{d}{dx}(\tilde{p}\tilde{r})^{-1/4}   $$
and $z(t)$ in \eqref{liouville} we obtain the following second order ODE
\begin{equation}
\label{normal}
\ddot{z}(t)+\left[Q_{1}(t)+R(t)  \right]z(t)=0,\qquad (0<t<t_{1}),
\end{equation}
with $Q_{1}(t)=\la-\frac{\tilde{q}(x)}{\tilde{r}(x)}$ and $t_{1}$, the value of $t$ at $x=L$.
On the other hand, from \eqref{liouville} we get
\begin{equation}
\label{boundary}
\begin{array}{rcl}
G(0)\tilde{p}(0)\dot{z}(0)&=&-\left[  \tilde{\al}_{1}-G'(0)   \right]\frac{1}{G(0)}z(0)\\
G(L)\tilde{p}(L)\dot{z}(t_{1})&=&-\left[  \tilde{\beta}_{1}-G'(L)   \right]\frac{1}{G(L)}z(t_{1})
\end{array}
\end{equation}
Finally, as $G(0)\tilde{p}(0)$ and $G(L)\tilde{p}(L)$ are positive constants, using the classical results for regular SLP given in \cite{eastham1970theory} Theorem 5.5.1, it follows the results in Theorem \eqref{asintotico} for the asymptotic behaviors of the eigenvalues of \eqref{problema}.
$\square$
\clearpage

\section{Applications. Studies on stability in three-layer Hele-Shaw flows}
\label{applications}
Using the results presented in the previous sections, in this part we consider a linear stability problem of the interfaces of a three-layer Hele-Shaw flow. This problem is a model to study a secondary oil recovery process and was presented in  \cite{gorell1983theory}.
This class of problems arise when oil is displaced by water through a porous medium producing the "fingering" phenomenon.
Only for the sake of completeness, we give an introduction to the deduction of the model.
See \cite{carasso1998optimal} for more details on the deduction an physical considerations.

In \cite{gorell1983theory} the authors S. Gorell and M. Homsy consider the displacement process by less viscous fluid containing a solute, and present a policy under which the Saffman-Taylor-Chouke instability can be minimized.
As the concentration of the solute is not constant and knowing that is possible to relate the concentration to the  viscosity, the porous medium can be considered saturated by three immiscible fluid: water, polymer-solute and oil.

The equations which govern the flow through a porous medium are
\begin{equation}
\label{sys}
\nabla\cdot\overrightarrow{V}=0,\qquad \nabla P=-\mu \overrightarrow{V},\qquad \frac{\pp \mu}{\pp t}+\overrightarrow{V}\cdot \nabla \mu=0,
\end{equation}
i.e.
the conservation of mass, the Darcy's law and the advection of viscosity for solute under the assumption that the adsorption, dispersion and diffusion are neglected.
This system admits the following steady displacement solution
\begin{equation}
\label{basicsol}
u=U;\quad v=0;\quad \mu=\mu_{0}(x_{1}-Ut);\quad P=-U\int_{x_{0}}^{x_{1}}\mu_{0}(x'-Ut)dx'=P_{0}
\end{equation}
For the stability analysis of this solution, linear perturbations are considered, obtaining
\begin{equation}
\label{sysli}
\nabla\cdot\overrightarrow{V'}=0,\qquad \frac{\pp P'}{\pp x}=-\mu' U-\mu_{0}u',\qquad \frac{\pp P'}{\pp y}=-\mu_{0}v' ,\qquad \frac{\pp \mu'}{\pp t}+u'\frac{d\mu_{0}}{dx}=0,
\end{equation}
where the coordinate $x_{1}$ was transformed into the moving reference frame $x=x_{1}-Ut$  and $u',\ v',\ P'$ and $\mu_{0}'$ denotes the component of the linear perturbation.
In what follows we consider that the lines $x=0$ and $x=L$ represent the interface between the fluids.
Since (\ref{sysli}) is linear, the perturbations can be represented by its Fourier integral.

Considering a typical wave component of the form (\emph{normal mode})
\begin{equation}
\label{ansatz}
 (u',v',p',\mu')=(f(x),\tau(x),m(x),n(x))\cdot e^{iky+\sigma t},
\end{equation}
where $k\in\rr$ denotes the wave number and $\sigma\in\cc$ denotes the growth rate of the perturbation of our basic configuration, and assuming that they are sectionally smooth functions, the ansatz (\ref{ansatz}) is consistent with (\ref{sysli}) provided
\begin{equation}
\label{compatible}
\tau(s)=ik^{-1}f'(s),\quad m(s)=-k^{-2} \mu_{0}(x)f'(x),\quad n(x)=-\sigma^{-1}\mu'_{0}(x)f(x),
\end{equation}
with $f(x)$ satisfying the following ODE
$$(\mu_{0}(x) f')'-k^{2}\mu_{0}(x) f  =-\frac{k^{2}U\mu'_{0}(x)}{\sigma}f,\qquad \textrm{when}\ x\neq0,\ x\neq L  \textrm{and}\ \sigma\neq 0.$$
 Here $()'$ denote the derivative with respect $x$.

Considering the following notations
\begin{equation}
\label{coeff}
\begin{array}{ccl}
\alpha_{1}(k)&=&k^{2}\left(\frac{S}{U}k^{2}+\mu_{1}-\mu_{0}(0^{+})   \right)\\
\alpha_{2}(k)&=&\mu_{1}k\\
\beta_{1}(k) &=&k^{2}\left(\mu_{2}-\mu_{0}(L^{-})-\frac{T}{U}k^{2}  \right)\\
\beta_{2}(2) &=&\mu_{2}k
\end{array}
\end{equation}
where $\mu_{1},\ \mu_{2}$ denotes the viscosity of water and oil respectively, $S$ the interfacial tension between water-polymer and $T$  between polymer-oil, the  dynamic and kinematic conditions to material points in the interfaces can be represented as follows
$$\mu_{0}(0)f'(0)=\left(\alpha_{1}(k)\frac{U}{\sigma} +\alpha_{2}(k)\right)f(0);\qquad  \mu_{0}(L)f'(L)=\left(\beta_{1}(k)\frac{U}{\sigma} -\beta_{2}(k)\right)f(L)$$
See Section 2 in \cite{gorell1983theory} for details on the above approximations.

Finally, considering
\begin{equation}
\label{growth}
\lambda=\frac{U}{\sigma},
\end{equation}
a spectral parameter, we obtain the following boundary value problem with the spectral parameter in the boundary condition
\begin{equation}
\label{problemintro}
(P)\left\lbrace\begin{array}{rcl}
 \dis{(\mu_{0}(x)f')'-(k^{2}\mu(x)-\lambda k^{2}\mu'_{0}(x)  )f}&=&0\qquad 0<x<L \\
                                          \mu_{0}(0)f'(0)&=&(\alpha_{1}(k)\lambda +\alpha_{2}(k))f(0)\\
                                          \mu_{0}(L)f'(L)&=&(\beta_{1}(k)\lambda -\beta_{2}(k))f(L)
            \end{array}\right.
\end{equation}
We remark that $\mu_{0}(x)$ is the viscosity of the polymer-solute and therefore the condition $\mu_{0}'(x)>0$ isn't a restrictive or unrealistic consideration.
Therefore, for problem \eqref{problemintro} we can use the results presented in this article.

This way, we can recognize the following behavior for the spectrum of (\ref{problemintro}):
\begin{coro}
\label{principalaplicacion}
The spectrum of (\ref{problemintro}) can be ordered in following cases:
\begin{itemize}
\item[i)] If $\alpha_{1}(k)<0$ and $\beta_{1}(k)>0$, then
 $$\quad 0<\lambda_{0}<\lambda_{1}<\lambda_{2}<\dots$$
\item[ii)] If $\alpha_{1}(k)\beta_{1}(k)\geq0$, then
$$\quad \lambda_{-0}<0<\lambda_{0}<\lambda_{1}<\lambda_{2}<\dots$$
\item[iii)]   If $\alpha_{1}(k)>0$ and $\beta_{1}(k)<0$, then
$$\quad\lambda_{-1}<\lambda_{-0}<0<\lambda_{0}<\lambda_{1}<\lambda_{2}<\dots$$
\end{itemize}
Moreover, the eigenfunction $f_{l}(x)$ has exactly $\vert l\vert$ zeros in $]0,L[$.
\end{coro}
The case i) is proved using the classical results presented in \cite{ince1962ordinary}; case ii) is presented in Appendix \ref{otroespectro} and finally, the case iii) is related with the work developed in this paper.

\begin{rem}
For the identification of the cases presented in \eqref{principalaplicacion}, de la definition of the physical parameters \eqref{coeff} we consider
\begin{equation}
\label{rangos}
\begin{array}{rcl}
\underline{k}^{2}&=&\min\left\lbrace  \frac{U}{S}(\mu_{0}(0^{+})-\mu_{1}),\     \frac{U}{T}(\mu_{2}-\mu_{0}(0^{+})) \right\rbrace\\
\\
\overline{k}^{2} &=&\max\left\lbrace  \frac{U}{S}(\mu_{0}(0^{+})-\mu_{1}),\     \frac{U}{T}(\mu_{2}-\mu_{0}(0^{+})) \right\rbrace
\end{array}
\end{equation}
Thus, the case i) is given for those $k<\underline{k}$, the case ii) for $k\in[\underline{k},\overline{k}]$, while iii) manifests for $k>\overline{k}$.
\end{rem}

Concerning the amplitude of the perturbative waves, from \eqref{growth}, when $k<\underline{k}$ we have that all growth rates are positive.
This might lead us to think that the most unstable case is for $k<\underline{k}$.
We remark that the degree of instability is set by the magnitude $\sigma_{0}=\frac{U}{\la_{0}}$.

In subsection \ref{subsec_numerico}, numerical computations are presented to provide an approximate understanding of the behavior of $\sigma_{0}$ when changing the wave number.
The aim of the numerical computations is to attain information on the dependency of the spectrum on the physical parameters.
\newpage

\subsection{Numerical experiments for linear middle viscous profile.
}\label{subsec_numerico}
In this part we present numerical approximations of the eigenvalues obtained for the particular case of a linear profile for the viscosity of the intermediate fluid, specifically, in \eqref{problemintro} we take $\mu_{0}(x)=ax+b$, where
$$ a=\frac{(\mu_{2}-\mu_{1})-(J_{1}+J_{2}) }{L};\qquad b=J_{1}+\frac{J_{2}+J_{1}-(\mu_{2}-\mu_{1})}{L}  $$
Here, we denote the viscosity jumps in the interfaces described by $x=0$ and $x=L$ as $J_{1}$ and $J_{2}$.

The results are presented in Tables \ref{cuadro1} and \ref{cuadro2}.
The numerical approximations are obtained using the nonlinear equation \eqref{laecuacion}.
The values for the physical parameters are the following:
$$S=1;\ T=1;\ L=0.1;\ \mu_{1}=1;\ \mu_{2}=2;\  J_{1}=0.1;\  J_{2}=0.1$$
In the experiment's developed we studied the spectra of \eqref{problemintro} for the cases: $U=1$ (Table \ref{cuadro1}) and $U=10$ (Table \ref{cuadro2}).

The values of the spectral parameter coefficients on the boundary conditions are the following:
\begin{equation}
\label{condicionesdeborde}
\begin{array}{lcll}
\al_{1}(k)=k^{2}\left[k^{2}-0.1 \right];& & \beta_{1}(k)=k^{2}\left[ 0.1-k^{2} \right]&\textrm{(case $U=1$)}\\
\\
\al_{1}(k)=0.1 k^{2}\left[k^{2}-0.01 \right];& & \beta_{1}(k)=0.1 k^{2}\left[ 0.01-k^{2} \right]&\textrm{(case $U=10$)}
\end{array}
\end{equation}
For the linear profile case, the solution for the ODE in \eqref{problemintro} can be represented through of Kummer and Tricomi confluent hypergeometric functions:
\begin{equation}
\label{hyper}
 \left\lbrace           \begin{array}{ccl}
                         \Phi(\al,1;z)&=& 1+\dis{\sum_{l=0}^{\infty}\frac{(\al)_{l}}{(l!)^{2}}z^{l}}\\
                            \\
                         \Psi(\al,1;z)&=&\dis{ \frac{1}{\Gamma(\al)}\left( \Phi(\al,1;z)\ln(z)+\sum_{l=0}^{\infty}\frac{(\al)_{l}}{(l!)^{2}}\frac{\Gamma'(\al+l)}{\Gamma(\al+l)}z^{l}    \right)   }
                \end{array}\right.
\end{equation}
Here $(\al)_{k}$ denotes the Pochhammer symbol and $\Gamma$ denotes the Euler-gamma function.
The fact that $h_{1}(\la)$ can be defined using these special function is direct from the consideration
\begin{equation}
\label{sus}
z=\frac{2k}{a}(ax+b);\qquad f(s)=e^{-z/2}g(z)
\end{equation}
Now, replacing in the ODE \eqref{problemintro}, the following Kummer equation is obtained
\begin{equation}
\label{kummer}
z\frac{d^{2}g}{dz^{2}}+(1-z)\frac{dg}{dz}-\left( \frac{1}{2}- \frac{k}{2}\la   \right)g=0,
\end{equation}
and therefore, the solution can be defined as follows
\begin{equation}
\label{gen}
g(z;\la)=C_{1}\Phi((1-\la k)/2,1;z)+C_{2}\Psi((1-\la k)/2,1;z)
\end{equation}
with $\Phi$ and $\Psi$ in \eqref{hyper} and $\al=(1-\la k)/2$.
\begin{table}[!hbt]
\begin{center}
\begin{tabular}[0.85\textwidth]{| c | c | c | c | c| c |c | c |}
  \hline
    $k$ &    $\alpha_{1}(k)$& $\beta_{1}(k)$  &$\la_{-1}$  & $\la_{-0}$   &$\la_{0}$     & $\la_{1}$      &  $\la_{2}$   \\
  \hline
  \hline
   1    &  0.9              & -0.9         & -6.46019      &  -3.27535    & 14.4968      & 68.1856         &  158.906          \\
   \hline
   2   &   15.6             &-15.6         & -0.51684      &   -0.29893   &  6.11938     & 19.6737         &   42.3671         \\
   \hline
   3   &   80.1             &  -80.1       &  -0.143895    &  -0.084287   &  3.81799     &  9.86499        &  19.9527         \\
   \hline
   4   &   254.4            & -254.4       & -0.0601436    & -0.0347156   & 2.96087      &  6.3777         & 12.0525         \\
   \hline
   5   &    622.5           &-622.5        & -0.0307741    & -0.0175412  & 2.55349       &   4.75521       &  8.38692        \\
   \hline
   6   &   1292.4           & -1292.4      & -0.017825     & -0.010068   & 2.32709       & 3.87244         & 6.39387          \\
   \hline
   7   &  2396.1            & -2396.1      & -0.0112375    & -0.00630497 & 2.18659       & 3.34028         & 5.19194       \\
   \hline
   8   &  4089.6            & -4089.6      &-0.007536      & -0.0042069  &  2.09175      &  2.99538        & 4.4122           \\
   \hline
   9   &  6552.9            & -6552.9      & -0.0052976    & -0.00294567 &  2.0233       & 2.75939         & 3.87815          \\
   \hline
   \hline
 \end{tabular}
 \caption{\label{cuadro1} Approximations considering $U=1$.
Here $S/U=1$ and $T/U=1$.
In our numerical experiment we use: Software \emph{Mathematica 9.0};  Library functions: Hypergeometric1F1[a,b,z]; HypergeometricU[a,b,z] for the  Kummer and Tricomi confluent hypergeometric functions respectively.}
\end{center}
\end{table}

\begin{table}[!hbt]
\begin{center}
\begin{tabular}[0.85\textwidth]{| c | c | c | c | c| c |c | c | c | c |c|}
  \hline
    $k$ &    $\alpha_{1}(k)$& $\beta_{1}(k)$  &$\la_{-1}$  & $\la_{-0}$   &$\la_{0}$     & $\la_{1}$      &  $\la_{2}$        \\
  \hline
  \hline
   1    &           0       &      0          &        -     &      -       &   4.96968    &26.7236         &        81.7087     \\
   \hline
   2   &  1.2               &   -1.2          &    -10.5131  & -6.14168     &  4.38068     &  16.2001       & 38.3316              \\
   \hline
   3   &   7.2              &    -7.2         &   -1.84609   &  -1.06909    & 3.51446      &  9.30549       &  19.3001             \\
   \hline
   4   &   24               &   -24           &  -0.677423   & -0.38866     &  2.88683     &  6.22874       &  11.8694            \\
   \hline
   5   &   60               &   -60           &  -0.329265   &  -0.187139   & 2.5302       &  4.70374       &  8.31946          \\
   \hline
   6   &    126             &  -126           &  -0.186085   &  -0.104951   & 2.31824      &  3.85155       &  6.36463             \\
   \hline
   7   &    235.2           &   -235.2       &   -0.115748   &  -0.0648894  &  2.18267     & 3.33682        & 5.17781               \\
   \hline
   8   &  403.2            &  -403.2         & -0.076998     &   -0.049627  &  2.08979     &  2.99072       &  4.40481               \\
   \hline
   9   &  648              &    -648         & -0.0538465    &  -0.0299318  &2.02221       &   2.75694      &   3.87404             \\
   \hline
   \hline
 \end{tabular}
 \caption{\label{cuadro2} Approximations considering $U=10$.
Here $S/U=0.1$ and $T/U=0.1$.
The approximation for the values of $\la_{l}$  }
\end{center}
\end{table}

In the results listing of Tables \ref{cuadro1} and \ref{cuadro2} we note a monotone behavior of the eigenvalues as the wave number increments.
We believe that this fact can be related to the asymptotic results of the eigenvalues as functions of the wave number $k$.
That, under the fact that we have $\beta_{1}(k)\sim -\frac{T}{U}k^{4}$ and $\beta_{2}(k)=\mu_{2}k$, and therefore have $h_{2}(\la,k)\sim -\frac{T}{U}k^{4}-\frac{\mu_{2}k}{\la}$.
This way, the graph intersection between $h_{1}(\la,k)$ and $h_{2}(\la,k)$ tend to $\mu_{n}(k)$, where  $\mu_{n}(k)$ is eigenvalue of the auxiliary problem \eqref{aux1}.

\section{Conclusions and comments}\label{sec_com}
Using classic tools on second order ODE ans elemental tools of mathematical analysis, we have obtained a set of results for the characterization of the eigenvalues and eigenfunctions of a SLP with spectral parameter in both boundary conditions.

For reaching our main results,we developed a list of lemmas that correspond to elemental technical adaptations in regular SLP, e.g. comparison  Sturm's theorem and separation theorem among other, see Section \ref{pre1}.
Concerning the techniques used for the description of the spectrum presented in other articles on the subject (see \cite{amara1999sturm,binding1994sturm} and their references), in our analysis we have considered a variant for the study of the characteristic equation, see equation \eqref{laecuacion}.
This variant corresponds to the definition of the function  $h_{1}(\la)$ in \eqref{funciones}.
The reason for this consideration is in the construction of the associated ODE through which we obtained the sign of the derivative and, therefore, the monotone behavior of $h_{1}(\la)$, see Equation \eqref{ODE} and Equation \eqref{sigder}.
This behavior of $h_{1}(\la)$ allowed us to order the real eigenvalues of \eqref{problema} bounding them as follows:
\begin{equation}
\label{orden2}
 \la_{-1}<\eta_{-0}<\la_{-0}<0<\la_{0}<\eta_{0}<\la_{1}<\eta_{1}<\la_{2}<\eta_{2}<\dots
\end{equation}
where $\eta_{l}$ are the eigenvalues of an auxiliary SLP given in  \eqref{aux1}.
About this auxiliary problem, we note that the fixed boundary condition in $x=0$ allows us to use directly some classical tools, e.g. Sturm's comparison criteria.

In order to have that the spectrum of our problem \eqref{problema} is a real subset, we use a Crum-Darboux type transformation, see Corollary \ref{soloreales}.
Also, through a Liouville transformation we have shown that the main problem \eqref{problema} has a regular SLP associated, and have obtained information on the spectrum through results on regular SLP, see Theorem \ref{asintotico}.

For the oscillatory results of the eigenfunctions, we have developed auxiliary results using the implicit function theorem.
On this results, we obtained a description of the behavior of the zeros of the functions that satisfy problem \eqref{main2}, see Lemmas \ref{zeros1} and \ref{zeros2}.
By means of this lemmas and the oscillatory results of Strum we proved the oscillatory behavior of the eigenfunctions, indicated in \ref{principal}, the separation Theorem \ref{principal2} (for the solutions of \eqref{main2}) and as a direct consequence of Corollary \ref{separacion}, which corresponds to the separation results for the eigenfunctions.

The results on \eqref{problema} have been applied in a three-layer Hele-Shaw flows model for the study of hydrodynamic stability of planar interfaces, see Section \ref{applications}.
This model was presented in \cite{gorell1983theory}, where the authors developed a theory which is able to describe the optimal policy which, if followed, minimizes the effects of the Saffman-Taylor-Chouke instability.

We comment that the initial motivation for our study on problem \eqref{problema} was to understand the non-regular SLP \eqref{problemintro}, therefore the hypotheses on the coefficient functions of the ODE and the signs of the coefficients on boundary conditions obey to that problem, see Corollary \ref{principalaplicacion}.
Concerning the hydrodynamic model \eqref{problemintro}, articles on the subject usually treat the case where $k<\underline{k}$, see \cite{carasso1998optimal} and references therein.
Out of the cases mentioned in Corollary \ref{problemintro} we have obtained a complete description of the behavior of the spectrum respective to the wave number $k$.
For an (introductory) analysis of the behavior of the growth rate of the perturbative waves, in Section \ref{subsec_numerico} we consider a particular case for $\mu_{0}(x)$.
Making numerical computations, we verified the results stated in our main results, synthesized in Corollary \eqref{problemintro}.
The list of approximations presented in Tables \ref{cuadro1} and \ref{cuadro2} is backed up using the oscillatory results on the eigenfunctions.
Otherwise we would have no certainty that the order of the eigenvalues is correct.

Generally, the instability and their description is related with the difference of viscosity between phases, see \cite{chuoke1959instability}, \cite{saffman1958penetration}.
In the numerical results we noted a high degree of dependency of the spectrum of \eqref{problemintro} and the physical parameters $L$, $S$, $T$, $U$, $\mu_{1}$ and $\mu_{2}$.
A more detailed analysis on this parameters requires other tools than the ones we have used in this work. We look to developing such tools in a future work.

Concerning the stability problem on a secondary oil recovery process, given the complexity of the phenomenon, we believe that the consideration of a linear profile for  $\mu_{0}(x)$ answers to elemental needs, and the  consideration of profiles $\mu_{0}(x)$ may be lacking of practical interest unless the explicit or numerical solutions allow a deeper analysis of the behavior of $\sigma_{n}(k)$ respective to the physical parameters mentioned.
To understand the degree of sensibility in the answer of the system to perturbations, a detailed analysis of the behavior of $\sigma_{0}=\sigma_{0}(k)$ must be developed.
We believe that by means of elemental analysis tools (e.g IFT) it's possible to obtain results for the description of $\sigma_{0}(k)$. We look to be able to develop such tools in a future article.

Concerning the geometrical approach to the spectrum and its relation with the oscillatory results of the eigenfunctions, we believe that it has reach for more general problems, and for problems that have been treated with other class of tools.
In this direction we put our tools to test 
That problem was addressed in the article \cite{daripa2008studies}, where the existence of a particular class of waves, called neutral waves is stated.
In Appendix \ref{susec_neutral} we prove that such class of wave doesn't exist, and explore the question whether the curves $\sigma_{n}(k)$ can intersect.
Such phenomenon would be in violation of the oscillator results of the eigenfunction, see Table \ref{espectros1}.

Finally, we believe that our tools can be adapted for superior order non-selfadjoint SLP, like the ones emerging on problems of hydrodynamics, see \cite{chandrasekhar1954characteristic, chandrasekhar1961hydrodynamic}.

\section*{\bf Acknowledgments:}\
The work of the second author (Oscar Orellana) was supported in part by Fondo Nacional de Desarrollo Cent\'ifico y Tecnologico (FONDECYT) under grant 1141260 and Universidad T\'ecnica Federico Santa Mar\'ia, Valpara\'iso, Chile.
The statements made herein are solely the responsibility of the authors.

\begin{appendices}
\label{apendice}
\section{On the spectrum for constant viscous middle profile.}\label{susec_neutral}

The oscillatory results of eigenfunctions can be used to determine if the spectrum of the eigenvalue problem is well defined. In this direction, The aim in this part is to use the geometrical approach and the oscillatory results on the eigenvalues and eigenfunctions to fully understand the spectrum of problem \eqref{problemintro} for the particular case $\mu_{0}(x)=\mu$ a constant function in $[0,L]$.\\
This model was studied in \cite{daripa2008studies} to obtain an upper bound for the growth rate in a simple unstable model of multi-layer Hele-Shaw flows, among other results.
.

For this particular case, \eqref{problemintro} is rewritten as follows:
\begin{equation}
\label{modelconstant}
\left\lbrace\begin{array}{rcl}
 f''- k^{2}f&=&0\qquad 0<x<L \\
     p(0)f'(0)&=&(\alpha_{1}(k)\lambda +\alpha_{2}(k))f(0)\\
       p(L)f'(L)&=&(\beta_{1}(k)\lambda -\beta_{2}(k))f(L)
            \end{array}\right.
\end{equation}
Since the ODE is non oscillatory, there are no eigenvalues $\la_{l}$ with $l\geq 2$.
This fact becomes clear by considering
\begin{equation}
\label{solmodelo}
f(x;\la)=\mu\cosh(kx)+\frac{(\al_{1}(k)\la+\al_{2}(k))}{k}\sinh(kx),
\end{equation}
the solution of the ODE in (\ref{modelconstant}), satisfying the boundary condition in $x=0$.
Thus
\begin{equation}
\label{fun1}
h_{1}(\la)=\frac{k}{\la}\frac{k\mu\sinh(kL)+(\al_{1}(k)\la+\al_{2}(k))\cosh(kL)}{k\mu\cosh(kL)+(\al_{1}(k)\la+\al_{2}(k))\sinh(kL)}.
\end{equation}
In \cite{daripa2008studies} the author presents the dispersion relation between the parameters $\la$ and $k$ y develops an analysis on the stability.
We focus mainly on the postulation on the existence of neutral waves (see \cite{daripa2008studies} III. A).

Using the functions in \ref{funciones}, the dispersion relation can be obtained using the geometric through the algebraic equation
$$h_{1}(\la)-h_{2}(\la)=0.$$
On the other hand, we note that
\begin{equation*}
\lim_{\la\to\infty}h_{1}(\la)=\lim_{\la\to-\infty}h_{1}(\la)=0.
\end{equation*}
Using the branches of $h_{1}(\la)$, in this part we present a complete characterization of the spectrum of the model problem (\ref{modelconstant}) and clarify some comments presented in \cite{daripa2008studies}.

Following similar arguments as the proofs in sections above, the next step corresponds to the search of the spectrum of the auxiliary problem.
To this end, when $\al_{1}(k)\neq 0$, the equation $f(0;\la)=0$ is equivalent
\begin{equation}
\label{valores}
\la=-\frac{1}{\al_{1}(k)}\left( \frac{\mu k}{\tanh(kL)}+\al_{2}(k)\right).
\end{equation}
As the bracket term is positive then, the sign of $\la$ is determined using the sign of boundary coefficient $\al_{1}(k)$.
Thus,
\begin{itemize}
\item If $\al_{1}(k)>0$ then $\la=\eta_{-0}$.
\item If $\al_{1}(k)<0$ then $\la=\eta_{0}$.
\end{itemize}
From the form of $f(s;\la)$ in (\ref{solmodelo}), for  cases $p(s)=\mu$ constant viscosity profile, the spectrum of the (sub) auxiliary problem
\begin{equation*}
\left\lbrace\begin{array}{rcl}
 f''- k^{2}f&=&0\qquad 0<x<L \\
     f(0)&=&0\\
     f(L)&=&0
            \end{array}\right.
\end{equation*}
 consists only of one element and this eigenvalue is determinate in (\ref{valores}) and  $h_{1}(\la)$ has only three branches.\\
In the Table \ref{espectros1} we present the spectrum of (\ref{modelconstant}) at the cases $\al_{1}(k)\neq 0$.
\begin{table}[!hbt]
\begin{center}
\begin{tabular}[0.85\textwidth]{ c  c  c  c  c  c  }
  \hline
  \hline
  $\al_{1}(k)$         &          $\beta_{1}(k)$       &  $\la_{-1}$  &    $\la_{-0}$  & $\la_{0}$ & $\la_{1}$     \\
  \hline

  \hline
  negative             & positive             &            not        &       not       & exist   & exist
  \\
  negative             &0                    &          not           &      not        & exist   & not
  \\
  negative             &negative             &          not           & exist           & exist   & not
  \\
   positive             & positive           &        not             &  exist          & exist   & not
  \\
  positive             &0                    &        not             & exist            & not    & not
  \\
  positive             &negative             &        exist           & exist            & not    & not

  \\

   \hline
   \hline
\end{tabular}
 \caption{\label{espectros1} Spectrum of (\ref{modelconstant}) at $\al_{1}(k)\neq 0$.}
\end{center}
\end{table}
We remark that the case $\al_{1}(k)=0$ is a critical case.
In this case, $f(s;\la)$ in (\ref{solmodelo}) and $f'(s;\la)$ are both positive functions and $h_{1}(\la)$ has only two branches, $\mathcal{B}_{-0}=]-\infty,0[$ and $\mathcal{B}_{0}=]0,\infty[$.
Thus, when $\al_{1}(k)=0 $ we get
\begin{itemize}
\item[] If $\beta_{1}(k)>0$ then $\Lambda=\lbrace \la_{0} \rbrace$
\item[] If $\beta_{1}(k)=0$ then $\Lambda=\emptyset$
\item[] If $\beta_{1}(k)<0$ then $\Lambda=\lbrace \la_{-0} \rbrace$
\end{itemize}
Finally, from Table \ref{espectros1} we note that the neutral waves postulated in \cite{daripa2008studies} do not exist.

\section{Attainment of the regular SLP \eqref{problem2}}\label{regularslp}
Through the following steps we obtain a  regular SLP associated with $g(x)$.
Here we use the ideas given in \cite{binding2004transformation}.
For the sake of completeness, we present the steps for the deduction of \eqref{problem2}.

From \eqref{transformacion}, we have that
\begin{equation}\label{derivada}
\begin{array}{rcl}
\displaystyle{\frac{dg}{dx}}&=&\displaystyle{(pf')'-pf'\frac{y_{0}'}{y_{0}}-f(py'_{0}/y_{0})'}\\
             &=& \displaystyle{(q-\la r)f-f' p\frac{y_{0}'}{y_{0}}-f\left\lbrace \frac{(py'_{0})'}{y_{0}}-p\left(\frac{y'_{0}}{y_{0}}\right)^{2} \right\rbrace} \\
             &=& \displaystyle{(q-\la r)f-  f' p\frac{y_{0}'}{y_{0}}-f\left\lbrace  (q-\la_{0}r)- p\left(\frac{y'_{0}}{y_{0}}\right)^{2}    \right\rbrace}\\
             &=& \displaystyle{(\la_{0}-\la)r f-\frac{y'_{0}}{y_{0}}\left\lbrace  p f'-pf\frac{y'_{0}}{y_{0}}  \right\rbrace}\\
             &=& \displaystyle{ (\la_{0}-\la)r f-\frac{y'_{0}}{y_{0}}g}.
\end{array}
\end{equation}
Dividing by $r$ and taking the derivative, it follows
\begin{equation}\label{derivada2}
\begin{array}{rcl}
\displaystyle{\frac{d}{dx}\left(\frac{1}{r}\frac{dg}{dx}\right)}&=&\displaystyle{(\la_{0}-\la)f'-\left(\frac{1}{r}\frac{y'_{0}}{y_{0}}\right)' g-
\frac{1}{r}\frac{y'_{0}}{y_{0}}\frac{dg}{dx}}\\
            &=&\displaystyle{ (\la_{0}-\la)f'-\left(\frac{1}{r}\frac{y'_{0}}{y_{0}}\right)' g-
            \frac{1}{r}\frac{y'_{0}}{y_{0}}\left\lbrace (\la-\la_{0}) r f-\frac{y'_{0}}{y_{0}}g  \right\rbrace }\\
            &=&\displaystyle{\left\lbrace  (\la_{0}-\la)\left(f' - \frac{y'_{0}}{y_{0}}f\right)\right\rbrace +\left\lbrace -\left(\frac{1}{r}\frac{y'_{0}}{y_{0}}\right)' +\frac{1}{r}\left(\frac{y'_{0}}{y_{0}}\right)^{2}
              \right\rbrace}g\\
              &=&\displaystyle{\left\lbrace  (\la_{0}-\la)\frac{p(x)}{p(x)}\left(f' - \frac{y'_{0}}{y_{0}}f\right)\right\rbrace +\left\lbrace -\left(\frac{1}{r}\frac{y'_{0}}{y_{0}}\right)' +\frac{1}{r}\left(\frac{y'_{0}}{y_{0}}\right)^{2}
              \right\rbrace}g\\
            &=&\displaystyle{\left\lbrace  \frac{1}{p(x)}(\la_{0}-\la)\left(p(x)f' - \frac{y'_{0}}{y_{0}}\mu_{0}(x)f\right)\right\rbrace +\left\lbrace -\left(\frac{1}{r}\frac{y'_{0}}{y_{0}}\right)' +\frac{1}{r}\left(\frac{y'_{0}}{y_{0}}\right)^{2}
              \right\rbrace}g\\
            &=&\displaystyle{\left\lbrace  \frac{\la_{0}}{p(x)} -\left(\frac{1}{r}\frac{y'_{0}}{y_{0}}\right)' +\frac{1}{r}\left(\frac{y'_{0}}{y_{0}}\right)^{2}-\frac{1}{p(x)}\la
              \right\rbrace}g
\end{array}
\end{equation}

\newpage
Now, we obtain the boundary condition  for $g(x)$.
At $x=0$ we have
\begin{equation}\label{condicion1}
\begin{array}{rcl}
g'(0)&=&\displaystyle{(\la_{0}-\la)r(0)f(0)-\frac{y'_{0}(0)}{y_{0}(0)}g(0)}\\
     &=&\displaystyle{r(0)\left\lbrace  \la_{0}f(0)-\la f(0;\la)      \right\rbrace-\frac{y'_{0}(0)}{y_{0}(0)}g(0)}\\
     &=&\displaystyle{r(0)\left\lbrace \la_{0}f(0)-\frac{\al_{1}\la f(0)}{\al_{1}}   \right\rbrace-\frac{y'_{0}(0)}{y_{0}(0)}g(0)}\\
     &=&\displaystyle{\frac{r(0)}{\al_{1}}\left\lbrace \al_{1}\la_{0}f(0)-(\al_{1}\la f(0))   \right\rbrace-\frac{y'_{0}(0)}{y_{0}(0)}g(0)}\\
     &=&\displaystyle{\frac{r(0)}{\al_{1}}\left\lbrace \la_{0}\al_{1}f(0)-(p(0)f'(0)-\al_{2}f(0) )   \right\rbrace-
     \frac{y'_{0}(0)}{y_{0}(0)}g(0)}\\
     &=&\displaystyle{\frac{r(0)}{\al_{1}}\left\lbrace \frac{(\la_{0}\al_{1}y_{0}(0))}{y_{0}(0)}
     f(0)-(p(0)f'(0)-\al_{2}f(0) )   \right\rbrace- \frac{y'_{0}(0)}{y_{0}(0)}g(0)}\\
     &=&\displaystyle{\frac{r(0)}{\al_{1}}\left\lbrace \frac{p(0)y'_{0}(0)-\al_{2}y_{0}(0)}{y_{0}(0)}
     f(0)-(p(0)f'(0)-\al_{2}f(0) )   \right\rbrace- \frac{y'_{0}(0)}{y_{0}(0)}g(0)}\\
     &=&\displaystyle{\frac{r(0)}{\al_{1}}\left\lbrace \frac{p(0)y'_{0}(0)}{y_{0}(0)}f(0)-\al_{2}f(0)-
     (p(0)f'(0)-\al_{2}f(0) )   \right\rbrace- \frac{y'_{0}(0)}{y_{0}(0)}g(0)}\\
      &=&\displaystyle{\frac{r(0)}{\al_{1}}\left\lbrace \frac{p(0)y'_{0}(0)}{y_{0}(0)}f(0)-
     p(0)f'(0)   \right\rbrace- \frac{y'_{0}(0)}{y_{0}(0)}g(0)}\\
     &=&\displaystyle{\frac{r(0)}{\al_{1}}\left\lbrace -\left(-\frac{p(0)y'_{0}(0)}{y_{0}(0)}f(0)+
     p(0)f'(0)\right)   \right\rbrace- \frac{y'_{0}(0)}{y_{0}(0)}g(0)}\\
     &=&\displaystyle{-\frac{r(0)}{\al_{1}}g(0)- \frac{y'_{0}(0)}{y_{0}(0)}g(0)}\\
     &=&\displaystyle{-\left\lbrace\frac{r(0)}{\al_{1}}+ \frac{y'_{0}(0)}{y_{0}(0)}\right\rbrace g(0)}.
\end{array}
\end{equation}
Similar to $x=L$, we get
\begin{equation}
\label{condicion2}
g'(L)=-\left\lbrace \frac{r(L)}{\beta_{1}}+\frac{y'_{0}(L)}{y_{0}(L)}  \right\rbrace g(L)
\end{equation}

\section{Proof of Theorem \ref{secundario}}
\label{demoaux}
We remark that  \eqref{aux1} has a fixed boundary condition, and therefore some elemental result are directly applicable.

In order to obtain a direct similitude with \eqref{problema}, we consider the following transformations
\begin{equation}
\label{cambio}
z=L-x;\qquad \tilde{y}(z)= y(x),
\end{equation}
obtaining the following non-regular SLP equivalent
\begin{equation}
\label{aux1eq}
\left\lbrace\begin{array}{rcl}
 \dis{(\tilde{p}(z)y')'-(\tilde{q}(z)- \lambda \tilde{r}(z) )\tilde{y}}&=&0\qquad 0<z<L \\
                                                            \tilde{y}(0)&=&0\\
                                                            \tilde{ y}'(L)&=&(-\alpha_{1}\lambda -\alpha_{2})\tilde{y}(L)\\
            \end{array}\right.
\end{equation}
Under the condition that $\al_{2}>0$, we have that $\tilde{h}_{2}(\la)=-\al_{1}-\frac{\al_{2}}{\la}$ is strictly increasing on each branch $]-\infty,0[$ and $]0,\infty[$.

Similar to our analysis for the study of the eigenvalues of \eqref{problema}, for $\tilde{y}(z;\la)$, solution of the ODE in \eqref{aux1eq} satisfying the boundary condition $\tilde{y}(0)=0$, we consider the function
\begin{equation}
\label{funcion2}
\tilde{h}_{1}(\la)=\tilde{p}(L)\frac{\tilde{y}'(L;\la)}{\la\tilde{y}(L;\la)}
\end{equation}
defined for all $\la\neq 0$ such that $\tilde{y}(L;\la)\neq 0$. This way, we consider the following auxiliary regular SLP for \eqref{aux1eq}
 \begin{equation}
\label{aux2eq}
\left\lbrace\begin{array}{rcl}
 \dis{(\tilde{p}(z)y')'-(\tilde{q}(z)- \lambda \tilde{r}(z) )\tilde{y}}&=&0\qquad 0<z<L \\
                                                            \tilde{y}(0)&=&0\\
                                                            \tilde{ y}(L)&=&0\\
            \end{array}\right.
\end{equation}
Using classical results (see Theorem II Ince \cite{ince1962ordinary} 10.6), we know that the spectrum of \eqref{aux2eq} can be ordered as follows:
$$ \tilde{\eta}_{0}<\tilde{\eta}_{1}<\tilde{\eta}_{3}<\dots \nearrow \infty,$$
where the corresponding eigenfunction $\tilde{y}_{n}=\tilde{y}(z;\tilde{\eta}_{n})$ has exactly $n$ zeros in $]0,L[$. Moreover, concerning the boundary conditions on \eqref{aux1eq}, for $\tilde{y}_{0}$ we know there exists certain $z^{*}$ such that the function $\tilde{y}_{0}$ reaches a maximum in $z^{*}$ (if $\tilde{y}_{0}<0$ we take the minimum and reach the same result). This way:
$$\tilde{p}(z^{*})\tilde{y}''_{0}(z^{*})=(\tilde{q}(z^{*})-\tilde{\eta}_{0}\tilde{r}(z^{*}))\tilde{y}_{0}(z^{*})=0$$
Knowing that the coefficients are positive, we get that $\tilde{\eta}_{0}>0$. Then, we have the following branches for the definitions of function $\tilde{h}_{1}(\la)$:
$$ \tilde{B}_{-0}=]-\infty,0[;\ \tilde{B}_{0}=]0,\tilde{\eta}_{0}[;\ \tilde{B}_{1}=]\tilde{\eta}_{1},\tilde{\eta}_{2}[;\dots  $$
Knowing $\tilde{y}(0;\la)=0$, through Sturm's second comparison theorem, for each $\tilde{B}_{n}$ ($n=0,1,2,\dots$) we get that $\tilde{h}_{1}(\la)$  is a decreasing function where
$$ \lim_{\la\nearrow \tilde{\eta}_{n}}\tilde{h}_{1}(\la)=-\infty;\quad \lim_{\la\searrow \tilde{\eta}_{n}}\tilde{h}_{1}(\la)=\infty  $$
Thus, for each $\tilde{B}_{n}$ ($n=0,1,2,\dots$) there exists a unique $\eta_{0}$ solution of $\tilde{h}_{1}(\la)=\tilde{h}_{2}(\la)$. As $y(z;\la)$ in \eqref{funcion2} satisfies the ODE in \eqref{aux2eq} and the first boundary condition at $z=0$, from the definition of $\tilde{h}_{2}(\la)$ we get that $\eta_{n}$ is an eigenvalue of \eqref{aux1}. For the oscillatory results we use Sturm's first comparison theorem. Finally, using minorant and majorant Sturm problems, similar as the proof of Lemma \ref{sobrelafuncion}, we obtain the existence for $\eta_{-0}$  and the oscillatory result is obtained using the maximum principle. $\square$

\section{Spectrum for the case $\al_{1}\beta_{1}\geq 0$}\label{otroespectro}
In this part, we consider $\al_{1}$ and $\beta_{1}$ in \eqref{problema}, such that $\al_{1}\beta_{1}\geq 0$. We obtain that the spectrum in this case can be ordered as follows
\begin{equation}
\label{espectros2}
\la_{-0}<0<\la_{0}<\la_{1}<\dots\nearrow \infty
\end{equation}
We note that in this case, we can use Sturm's comparison theorem.

Without loss of generality we may consider $\al_{1}\leq 0$, otherwise, we may introduce a change of coordinates, like in \eqref{cambio}.
Under the consideration that $\al_{1}\leq 0$, the spectrum of the auxiliary problem \eqref{aux1} is ordered al follows:
$$ 0<\eta_{0}<\eta_{1}<\eta_{2}<\dots\nearrow \infty $$ 
Let $f(x;\la)$ be the solution of the initial value problem \eqref{main2} and let $h_{1}(\la)$ be defined in \eqref{funciones}. Since $\al_{1}<0$, for $\underline{\la}<\overline{\la}$, both in $\mathcal{B}_{n}$ (with $n=1,2,\dots$) from Sturm's first comparison theorem have that $f(x;\underline{\la})$ and $f(x;\overline{\la})$ have exactly $n$ zeros in $]0,L[$. Now, from Sturm's second comparison theorem we obtain the monotonic behavior of $h_{1}(\la)$. The existence of $\la_{0}<\la_{1}<\la_{2}<\dots\nearrow\infty$ is obtained in similar way to the proofs already presented. For the proof of the existence of $\la_{-0}$, similar arguments to the ones presented in Lemma \ref{sobrelafuncion} can be used.

\section{Relacion entre los espectros}\label{regularSLP}
In this appendix, we show that the spectra of problems \eqref{problema} and \eqref{problem2} are almost directly related.

We will denote the eigenvalues of \eqref{problem2} as follows
\begin{equation}
\label{espectroslpregular}
\tilde{\la}_{0}<\tilde{\la}_{1}<\tilde{\la}_{2}<\dots\nearrow \infty,
\end{equation}
with their respective eigenfunctions, $g_{m}(x)$.

We define
\begin{equation}
\label{funpropia}
w_{m}(x)=y_{0}(x)\int_{0}^{x}\frac{1}{y_{0}p}g_{m}ds +Cy_{0}(x),
\end{equation}
where $C$ is a constant to be determined. We have that
\begin{equation}
\label{relacion}
g_{m}(x)=p(x)w'_{m}-pw_{m}\frac{y'_{0}}{y_{0}}
\end{equation}
On the other hand, from the definition, it holds that
$$ p\frac{d w_{m}}{dx}=py'_{0}\left\lbrace  C+ \int_{0}^{x}\frac{1}{y_{0}p}g_{m}ds \right\rbrace +g_{m}.$$
Differentiating the previous equation we have
\begin{equation}
\label{step0}
\begin{array}{rcl}
\dis{\frac{d}{dx}(pw_{m}')}&=&\dis{(q-\la_{0}r)y_{0}\left\lbrace  C+ \int_{0}^{x}\frac{1}{y_{0}p}g_{m}ds \right\rbrace +\frac{y_{0}'}{y_{0}}g_{m}+\frac{d g_{m}}{dx}}\\
                     &=&\dis{(q-\la_{0}r) w_{m}+\frac{y_{0}'}{y_{0}}g_{m}+\frac{d g_{m}}{dx}}
\end{array}
\end{equation}
On the other hand, integrating the ODE \eqref{problem2} between 0 an $x$, we get
 \begin{equation}
 \label{step1}
 \begin{array}{rcl}
 \frac{1}{r(x)}\left\lbrace \frac{y_{0}'}{y_{0}}g_{m}+g'_{m} \right\rbrace&=& \frac{1}{r(0)}\left\lbrace \frac{y_{0}'(0)}{y_{0}(0)}g_{m}(0)+g'_{m}(0) \right\rbrace +\dis{(\la_{0}-\tilde{\la}_{m})\int_{0}^{x}\frac{g_{m}}{p} ds+  }\\
 && +\dis{\int_{0}^{x}\frac{1}{r}\frac{y'_{0}}{y_{0}}\left\lbrace \frac{y_{0}'}{y_{0}}g_{m}+g'_{m}  \right\rbrace ds}
 \end{array}
 \end{equation}
Now, using \eqref{relacion} we get
\begin{equation}
\label{step2}
\begin{array}{rcl}
\dis{  \int_{0}^{x}\frac{g_{m}}{p} ds }&=&\dis{  \int_{0}^{x}w'_{m}-w_{m} \frac{y'_{0}}{y_{0}}ds }\\
                                       &=&\dis{ w_{m}(x)-w_{m}(0)-  \int_{0}^{x}w_{m} \frac{y'_{0}}{y_{0}}ds      }
\end{array}
\end{equation}
Taking
\begin{equation}
\label{lafuncion}
F(x)=-\la_{0}rw_{m}+\frac{y_{0}'}{y_{0}}g_{m}+g'_{m},
\end{equation}
when we substitute \eqref{step2} in \eqref{step1}, we get
\begin{equation}
\label{step3}
F(x)=r(x)\left\lbrace \frac{1}{r(0)}F(0)+\tilde{\la}_{m} \right\rbrace+r(x)\left\lbrace \int_{0}^{x} \frac{y_{0}'}{y_{0}r}(F(s)+\tilde{\la}_{m}w_{m}(s))ds\right\rbrace-\tilde{\la}_{m}rw_{m}
\end{equation}
We take $w_{m}(0)$ (and therefore, some $C$ in \eqref{funpropia}) such that
\begin{equation}
\label{condicion0}
(\la_{0}-\tilde{\la}_{m})w_{m}(0)=\frac{1}{r(0)}g'_{m}(0)+\frac{1}{r(0)}\frac{y'_{0}(0)}{y_{0}(0)}
\end{equation}
Now, from \eqref{step0} and \eqref{step3} we have
\begin{equation}
\label{step4}
\dis{\frac{d}{dx}(pw_{m}')}-(q-r\tilde{\la}_{m})w_{m}=r(x)\left\lbrace \int_{0}^{x}\frac{y'_{0}}{r y_{0}}  (F(s)+\tilde{\la}_{m} w_{m}(s))ds  \right\rbrace
\end{equation}
Taking
\begin{equation}
\label{step5}
\tilde{F}(x)=\dis{\frac{d}{dx}(pw_{m}')}-(q-r\tilde{\la}_{m})w_{m}
\end{equation}
and the definitions in \eqref{lafuncion}, from \eqref{step4} we have
$$ \tilde{F}(x)=r(x)\dis{ \int_{0}^{x}\frac{y'_{0}}{p y_{0}} \tilde{F} ds}$$
Using the consideration \eqref{condicion0} we have that $\tilde{F}(0)=0$. Now we have the following initial vaule problem 
$$ \frac{d\tilde{F}}{dx}=\left\lbrace \frac{r'(x)}{r(x)}+\frac{y'_{0}(x)}{y_{0}(x)} \right\rbrace\tilde{F}(x);\qquad \tilde{F}(0)=0,  $$
and therefore, it holds that $\tilde{F}=0$. Then we have
\begin{equation}
\label{step6}
\dis{\frac{d}{dx}(pw_{m}')}-(q-r\tilde{\la}_{m})w_{m}=0
\end{equation}
Finally, knowing that it holds that $p(0)\frac{y'_{0}(0)}{y_{0}}=\la_{0}\al_{1}+\al_{2}$ and $p(L)\frac{y'_{0}(L)}{y_{0}(L)}=\beta_{1}\la-\beta_{2}$, from \eqref{condicionesregular} we have that
\begin{equation}
\label{step7}
p(0)w'_{m}(0)=(\al_{1}\tilde{\la}_{m}+\al_{2})w_{m}(0);\qquad p(L)w'_{m}(L)=(\beta_{1}\tilde{\la}_{m}-\beta_{2})w_{m}(L)
\end{equation}
This way, considering \eqref{step6} and \eqref{step7}, we have that $w_{m}(x)$ in \eqref{funpropia} is an eigenfunction of \eqref{problema}, with eigenvalue $\tilde{\la}_{m}$.
\end{appendices}

\bibliographystyle{plain}
\bibliography{bib}

\end{document}